\documentclass[11pt]{article}
\pdfoutput=1
\usepackage{lineno}
\usepackage{lmodern}
\usepackage{amsmath,amssymb,amsthm}
\usepackage{xspace}
\usepackage{microtype}
\usepackage{tikz}

\usepackage{pgfplots}
\usetikzlibrary{calc,positioning}
\pgfplotsset{compat=1.13}
\usepgfplotslibrary{fillbetween}

\usepackage[backgroundcolor=yellow]{todonotes}		
\usepackage{mathtools}
\usepackage[style=alphabetic,maxbibnames=100,maxcitenames=10,backend=bibtex]{biblatex}
\renewbibmacro{in:}{}

\usepackage[colorlinks = true,
			citecolor = blue,
			linkcolor = blue]{hyperref}

\newcommand{\dd}{\mathop{}\!\mathrm{d}}
\let\del\partial

\newcommand{\xaddspace}[0]{\mathchoice{\hspace{-0.8em}}
                {\hspace{-0.6em}}
                {\hspace{-0.4em}}
                {\hspace{-0.3em}}
                 }

\newcommand{\xint}[1]{\int\foreach \i in {2,...,#1}{\xaddspace\int}}
\newcommand{\myiint}[0]{\xint{2}}
\renewcommand{\iint}{\myiint}

\newcommand{\RR}{\mathring R}
\newcommand{\RRR}{\mathring{\bar R}}
\newcommand{\ootimes}{\mathbin{\mathring{\otimes}}}

\let\div\relax
\DeclareMathOperator{\div}{div}
\DeclareMathOperator{\tr}{Tr}
\newcommand{\vin}[0]{v^{\textup{in}}}
\newcommand{\pex}[0]{\textsf{\textit{p}}}
\newcommand{\vex}[0]{\textsf{\textit{v}}}
\newcommand{\vv}[0]{\bar v}

\newcommand{\ppp}[0]{\bar{\bar p}}
\newcommand{\pp}[0]{\bar p}
\newcommand{\wpq}[1][q+1]{w^{\textup{(p)}}_{#1}}
\newcommand{\wcq}[1][q+1]{w^{\textup{(c)}}_{#1}}
\newcommand{\Rnash}{R_{q+1}^\textup{Nash}}
\newcommand{\Rtransport}{R_{q+1}^\textup{trans}}
\newcommand{\Rosc}{R_{q+1}^\textup{osc}}
\newcommand{\matd}[1]{D_{t,#1}}
\newcommand{\TT}[0]{\mathsf{T}}
\newcommand{\vinreg}[0]{\bar\beta}%

\def\dashint{\,\ThisStyle{\ensurestackMath{%
  \stackinset{c}{.2\LMpt}{c}{.5\LMpt}{\SavedStyle-}{\SavedStyle\phantom{\int}}}%
  \setbox0=\hbox{$\SavedStyle\int\,$}\kern-\wd0}\int}
\def\ddashint{\,\ThisStyle{\ensurestackMath{%
  \stackinset{c}{.2\LMpt}{c}{.5\LMpt+.2\LMex}{\SavedStyle-}{%
    \stackinset{c}{.2\LMpt}{c}{.5\LMpt-.2\LMex}{\SavedStyle-}{%
      \SavedStyle\phantom{\int}}}}\setbox0=\hbox{$\SavedStyle\int\,$}\kern-\wd0}\int}

\usepackage{scalerel,stackengine}
\stackMath
\newcommand\widecheck[1]{%
\savestack{\tmpbox}{\stretchto{%
  \scaleto{%
    \scalerel*[\widthof{\ensuremath{#1}}]{\kern-.6pt\bigwedge\kern-.6pt}%
    {\rule[-\textheight/2]{1ex}{\textheight}}
  }{\textheight}%
}{0.5ex}}%
\stackon[1pt]{#1}{\scalebox{-1}{\tmpbox}}%
}

\DeclareSymbolFont{bbold}{U}{bbold}{m}{n}
\DeclareSymbolFontAlphabet{\mathbbold}{bbold}
\newcommand{\ind}[0]{\mathbbold{1}}

\newcommand{\coloneq}{\mathrel{\mathop:}=}
\newcommand{\eqcolon}{=\mathrel{\mathop:}}
\newcommand{\ii}{\textup{i}}
\newcommand{\ee}{\textup{e}}
\DeclareMathOperator{\curl}{curl}

\newcommand{\varEPS}{0.2}
\newcommand{\varTAU}{0.2}
\newcommand{\varNi}{2}
\newcommand{\varN}{4}
\newcommand{\varYOFF}{-8*\varTAU }
\newcommand{\varETASTART}{\varNi*\varTAU}

\newcommand{\colorR}{blue!50!white}
\newcommand{\colorEta}{blue!20!white}
\newcommand{\colorEtaMinusOne}{red}

\newcommand{\crvi}[1]{
    \varTAU/3 - (1-2*\varEPS)*\varTAU/3 
    }
\newcommand{\crvii}[1]{
    2*\varTAU/3 + (1-2*\varEPS)*\varTAU/3 
    }

\newcommand{\crviii}[1]{
    \varTAU/3 - (1-\varEPS)*\varTAU/3 
    + 2*\varEPS *\varTAU *sin(deg(2 * pi * #1))/3
    }
\newcommand{\crviv}[1]{
    2*\varTAU/3 + (1-\varEPS)*\varTAU/3 
    + 2*\varEPS *\varTAU *sin(deg(2 * pi * #1))/3
    }
    
\newcommand{\etaminusone}[1]{
    (1 - (#1)*(#1) ) * (1 - (#1)*(#1))
}

\newtheorem{thm}{Theorem}[section]
\newtheorem{cor}[thm]{Corollary}
\newtheorem{lem}[thm]{Lemma}
\newtheorem{prop}[thm]{Proposition}
\theoremstyle{definition}
\newtheorem{defn}[thm]{Definition}

\theoremstyle{remark}

\numberwithin{equation}{section}
\newcommand{\Holder}{H\"older\xspace}
\newcommand{\Gronwall}{Gr\"onwall}
\newcommand{\Szekelyhidi}{Szek\'elyhidi}

\DeclareMathOperator{\supp}{supp}
\bibliography{convexbib}

\renewcommand{\epsilon}{\varepsilon}
\begin{document}
\title{Infinitely many non-conservative solutions for the three-dimensional Euler equations with arbitrary initial data in $C^{1/3-\epsilon}$}
\author{Calvin Khor\footnote{
Beijing Computational Science Research Center,
Beijing 100193, China,
 \texttt{calvin\_khor@hotmail.com}}
\and  Changxing Miao\footnote{Institute of Applied Physics and Computational Mathematics,
P.O. Box 8009, Beijing 100088, P. R. China,  \texttt{miao\_changxing@iapcm.ac.cn}}
\and Weikui Ye\footnote{Institute of Applied Physics and Computational Mathematics,
P.O. Box 8009, Beijing 100088, P. R. China,  \texttt{904817751@qq.com}}}
\date{}
\maketitle
\begin{abstract}
Let $0<\beta<\bar\beta<1/3$. We construct infinitely many distributional solutions in $C^{\beta}_{x,t}$  to the three-dimensional Euler equations that do not conserve the energy, for a given initial data in $C^{\vinreg}$. We also show that there is some limited control on the increase in the energy for  $t>1$.
\end{abstract}
\allowdisplaybreaks
\section{Introduction}
\label{s:intro}
This paper studies the Cauchy problem for the 3D incompressible Euler equations,
\begin{equation}
\left\{ \begin{alignedat}{-1}
\del_t v+(v\cdot\nabla) v  +\nabla p   &=  0,
 \\
  \nabla \cdot v  &= 0,
  \\
  v\big|_{t=0} &= \vin ,
\end{alignedat}\right.  \label{e:euler}
\end{equation}
considered on the spatially periodic domain $\mathbb T^3 \coloneq  (\mathbb R/2\pi \mathbb Z)^3$. For $x\in\mathbb  T^3$ and $t\ge 0$, $v(x,t)\in\mathbb R^3$ is the fluid velocity, and $p(x,t)\in\mathbb R$ is the pressure field. We are specifically interested in the scenario of low regularity initial data $\vin\in C^{\vinreg}$ for $\vinreg<1/3$, in relation to the following celebrated result:
 \begin{thm}[Onsager's Conjecture]  Consider weak solutions of \eqref{e:euler} in $C^\beta(\mathbb T^3\times [0,T])$.
 \begin{enumerate}
    \item \cite{zbMATH00722097} If $\beta>1/3$, then the energy $\int_{\mathbb T^3} |v|^2 \dd t$ is conserved in time.
    \item \cite{zbMATH06976275} If  $\beta<1/3$, then there exist solutions that do not conserve energy.
\end{enumerate}
 \end{thm}
 The work \cite{zbMATH06976275} proving the negative part of Onsager's conjecture is the culmination of many refinements over many papers over many years, starting with the constructions of Scheffer \cite{zbMATH00427755} and Shnirelman \cite{zbMATH02019780}, followed by the discovery  \cite{zbMATH05710190} of the deep connection between a counterintuitive result in geometry (the Nash--Kuiper theorem \cite{zbMATH03093832,b1aaf276a2034cae0b984600fa938a91,4628c42acd42fbf46589e4536b5bf8f2}), and turbulence. In this seminal paper of De Lellis and Sz\'ekelyhidi \cite{zbMATH05710190}, they proved the existence of distributional solutions to the Euler equations in $L^\infty$ with compact support in spacetime.

After nearly half a decade, the authors switched from a `soft analysis' approach to `hard analysis' and managed \cite{zbMATH06210493} to show the existence of continuous solutions. This led to a flood of papers with the above mentioned work of Isett as a crowning result. One should also mention the streamlined proof of \cite{zbMATH07038033} that also allows you to prescribe an arbitrary positive energy profile. Our account of the history is extremely abridged; we refer readers to the surveys \cite{zbMATH07207865,zbMATH06996008,zbMATH07078776} for more detailed information.

Roughly in parallel with the above developments, a version of Onsager's conjecture for the Cauchy problem was also studied. As is well known, for  sufficiently smooth initial data, there exists a unique smooth solution for short times. The question whether such a solution is in fact global is a classical open problem. The above results were proven by showing that there are  $L^\infty$ \cite{zbMATH05710190} and $C^\beta$ ($\beta<1/3$) \cite{zbMATH06976275} solutions with compact support in time, and hence in particular that (so long as we are considering distributional solutions) there are multiple solutions starting from the identically zero flow. In fact, one expects convex integration to produce a set of solutions that is Baire-typical due to the h-principle (see e.g. \cite{zbMATH06710292}).

It is natural to predict that the same is true for other initial data. Daneri \cite{zbMATH06312794}, later with  \Szekelyhidi{}  \cite{zbMATH06710292} and  Runa \cite{zbMATH07370998}  proved (among other things) that the set of initial data in $C^{\vinreg}$  without uniqueness of  $C^\beta$ solutions ($0<\beta<\vinreg<1/3$) is $L^2(\mathbb T^3)$-dense. Rosa and Haffter \cite{Rosa2021DimensionOT} also showed that any $C^\infty$ initial data gives rise to uncountably many solutions.

The main result of our paper is the following theorem, which generalises the very strong ill-posedness above  to  \emph{all} initial data in $C^{\vinreg}$:
\begin{thm}\label{t:main}
Let $\vin \in C^{\vinreg} (\mathbb T^3)$ be divergence-free. Fix $0\le \beta<\vinreg <1/3$, and $T\in(1,\infty]$. Then there exists infinitely many smooth functions $e:[1,T]\to(\|\vin\|_{L^2}^2,\infty)$ with a corresponding a weak solution $v$ to \eqref{e:euler} on $[0,T]$ with initial data $v\big|_{t=0}=\vin$, such that $v\in C^\beta_{x,t}$, and for all $t\in [1,T]$, $\int_{\mathbb T^3} |v(x,t)|^2 \dd x = e(t)$.
\end{thm}
We achieve this by adjusting the convex integration scheme to take an approximate initial data, which recovers any chosen initial data in the limit. This shows the strongest possible ill-posedness for an equation that has solutions.

We note however that there are significant differences between  Theorem \ref{t:main} and the earlier mentioned works  \cite{zbMATH06312794,zbMATH06710292,zbMATH07370998}. Recall that weak-strong uniqueness (see \cite{zbMATH06999781}) holds for the Euler equations: energy dissipative weak solutions (called `admissible solutions') with smooth initial data are classical. The above works of Daneri and collaborators show that weak-strong uniqueness does not hold for all $C^{\vinreg}$ data. In contrast, the solutions that we construct  must \emph{increase the energy}: we have in effect traded some control the energy in exchange for control on the initial data. This energy increase is also present in the construction of \cite{Rosa2021DimensionOT}.

We still manage some control  over the energy, which is better than the energy control at a single point achieved by  \cite{zbMATH07370998}, for instance. More precisely, we have the following quantitative version of Theorem \ref{t:main}:
\begin{thm}\label{t:quant}
For the solutions constructed in Theorem \ref{t:main} with chosen parameters $a,b,\beta,\vinreg$, there exist constants $c,C  > 0 $ depending on the parameters of Proposition \ref{p:main-prop} (in particular diverging to $+\infty$ as $a\to\infty$, $b\to1$, or $\beta\to\vinreg$)  such that for each $\vin\in C^{\vinreg}$, the energy profile $e$ for times $t\ge 1$ can be chosen to be an arbitrary smooth function satisfying
\[ c< e(t) - \int_{\mathbb T^3}   |\vin|^2 \dd x < C \text{ for all $t\in[1,T]$.} \label{energy control}\]
\end{thm}
We manage this control of energy by transitioning at some arbitrary positive time (we have chosen $t=1$ for simplicity) from the construction of \cite{zbMATH06976275} (which allows velocity perturbations supported away from the initial time) to the squiggling cutoffs of \cite{zbMATH07038033}. 

It would be very interesting to know if each non-smooth initial data in $C^{\vinreg} $ has energy dissipative solutions, and to what extent can the energy be controlled, but these questions are beyond the scope of the current paper.


We finish this introductory section with a non-exhaustive list of papers where convex integration has  been adapted to show non-uniqueness for other equations, even equations with dissipation. For instance, there are results for Euler on $\mathbb R^3$ \cite{zbMATH06591765}, MHD \cite{zbMATH07327027,8c8b8100f398e1b6e7c50f216ac7247a},  active scalars \cite{zbMATH06919574}, 2D Euler \cite{c52c7ef0df543843a5767f09460b9580}, SQG \cite{zbMATH07360846}, Navier--Stokes  \cite{zbMATH07003146,buckmaster2020wild,zbMATH07201159}, transport equations \cite{zbMATH07370992}, stationary Navier--Stokes \cite{zbMATH07061535}, Boussinesq \cite{zbMATH07183858,zbMATH06766517,zbMATH06873447}.  Recently, Colombo \cite{brue2021nonuniqueness} showed the non-uniqueness of solutions to 2D Euler with vorticity in the Lorentz space $L^{1,p}$, $p>2$.  We  also mention  the paper \cite{dellis2020nonuniqueness} showing the non-uniqueness of globally dissipative solutions, which additionally solve the local energy inequality and strictly dissipate the energy (but at the moment, they are unable to achieve \Holder regularity beyond $C^{1/7}$).

The remainder of this paper is organised as follows. In Section \ref{s:notation}, we fix some notation. In Section \ref{s:outline}, we give the main iterative proposition, sketch its proof, and show how Theorems \ref{t:main} and \ref{t:quant} follow from it. Section \ref{s:construct} carries out the convex integration scheme (in particular that \eqref{e:velocity-diff} holds at step $q+1$ in Corollary \ref{c:estimates-for-wpq-wcq}). In Section \ref{s:est-stress}, we explain how the stress error terms are controlled, and in Section \ref{s:energy-iteration}, we show that we have the claimed control on the energy. Finally we collect some classical estimates on \Holder spaces in the Appendices.

\section{Notation}\label{s:notation}
For $\alpha\in(0,1)$, and $N\in\mathbb Z_{\ge0}$, $C^{N+\alpha}(X)$ denotes the usual \Holder{} space which we have defined in Appendix \ref{s:holder}. We will write $C^{N+\alpha}$ as shorthand for $C^{N+\alpha}(\mathbb T^3)$. For functions in $L^\infty(0,T; C^{N+\alpha})\eqcolon L^\infty C^{N+\alpha}$, we define the norm
 \[ \|f\|_{N+\alpha}\coloneq \|f\|_{L^\infty C^{N+\alpha}} .\]
 We also write $[f]_\alpha=\sup_{t\le T} [f(t,\cdot)]_{C^\alpha}$ for the \Holder seminorm.
Some classical estimates  for \Holder functions are collected in Appendix \ref{s:holder}.

We fix for the remainder of the paper an even non-negative bump function $\widetilde\phi\in C^\infty_c([-1,1];[0,\infty))$, i.e. $\tilde\phi(-t)=\tilde \phi(t)$, and define for each $\epsilon>0$ two sequences of  mollifiers with integral 1:
\begin{align}
     \phi(t) 
        &\coloneq \frac{\widetilde\phi(t)}{\int_{[-1,1]}\widetilde\phi (\tau)\dd \tau },
     &\phi_{\epsilon}(t) 
        &\coloneq \frac1{\epsilon} \phi\left(\frac t\epsilon\right)\label{e:defn-mollifier-t},
     \\
    \psi(x)
        &\coloneq\frac{\widetilde\phi(|x|)}{\int_{\mathbb B_1(0)}\widetilde\phi (|y|)\dd y},
    &\psi_\epsilon(x) 
            &\coloneq \frac1{\epsilon^3} \psi\left(\frac{x}\epsilon\right). \label{e:defn-mollifier-x}
\end{align}

We use the following distributional notion of a weak solution:
\begin{defn} Let $T\ge0$ and let $\vin\in C^{\vinreg}$ for some $\vinreg>0$ be divergence-free (in the sense of distributions) and have zero mean, i.e. $\int_{\mathbb T^3} \vin \dd x=0$. We say that $v\in C^0(\mathbb T^3 \times [0,T]) $ is a \emph{weak solution on $[0,T)$ to \eqref{e:euler} with initial data $\vin$} if $v(\cdot,t)$ has zero mean, is divergence-free, and  for all divergence-free test functions $f\in C^\infty_c([0,T)\times \mathbb T^3;\mathbb R^3)$,    \begin{gather*}
    \int_0^T \int_{\mathbb T^3} v\cdot (\del_t+v\cdot\nabla)f \dd x \dd t =-\int_{\mathbb T^3} \vin(x) f(0,x)\dd x.
\end{gather*}
\end{defn}
The zero mean condition above is for convenience, and is easily removed via the Gallilean symmetry $\widetilde v(x,t)=v(x+t \bar v,t)-\bar v$ with $\bar v = \int_{\mathbb T^3} v \dd x\in\mathbb R^3 $. The pressure $p$ which does not appear in the weak formulation is recovered by solving the Poisson equation $-\Delta p = \div\div(v\otimes v)$.  Further background on fluid mechanics can be found for instance in \cite{zbMATH01644218}.

\section{Outline of the convex integration scheme}
\label{s:outline}
\subsection{Parameters and their restrictions}\label{ss:params}
For all $q\ge 0$ and given parameters $\alpha\in(0,1)$, $a\in(2,\infty)$, $b\in(1,2)$,  $0<\beta<\vinreg<1/3$,
we define ($\lceil x\rceil$ denotes the ceiling function)
\begin{alignat*}{10}
    \lambda_q &\coloneq 2\pi \left\lceil a^{b^q}\right\rceil > 1 ,&\qquad   \delta_q &\coloneq \lambda_q^{-2\beta} < 1.
\end{alignat*}
Note that $2\pi \le \frac{\lambda_q}{a^{b^q}} \le 4\pi $. We restrict $b<2$ so that
$ \frac1{50}\le \frac{\lambda_{q+1}}{\lambda_q^b} \le 2.$ We will require that $\alpha\ll1$  and $a\gg 1$ so that
\begin{align}
     2\alpha < \beta (b-1),\quad  a > 50^{\beta/\alpha }. \label{e:params0}
\end{align}
These imply that $\lambda_q^{3\alpha} \le \frac{\delta_q^{3/2}}{\delta_{q+1}^{3/2}}  \le  \frac{\lambda_{q+1}}{\lambda_q}$, as in \cite{zbMATH07038033}:
\[ 
\lambda_q^{3\alpha} \le \frac{\lambda_q^{6\alpha}}{50^{3\beta}} \le \frac{\lambda_q^{3\beta(b-1)}}{50^{3\beta}} \le \left(\frac{\lambda_{q+1}}{\lambda_{q}} \right)^{3\beta} = \frac{\delta_q^{3/2}}{\delta_{q+1}^{3/2}} \le \frac{\lambda_{q+1}}{\lambda_q}.\]
We also define the mollification parameter $\ell_q$ by
\begin{align}
    \ell_q \coloneq \frac{\delta_{q+1}^{1/2}}{\delta_q^{1/2} \lambda_q^{1+3\alpha/2}}  \in \Big(\frac12\lambda_q^{-1-\frac{3\alpha}2-(b-1)\beta},\lambda_q^{-1-\frac{3\alpha}2}\Big) \label{e:ell}.
\end{align}
The lower bound on $\ell_q$ follows from $\frac12\lambda_q^{1-b} \le \lambda_q\lambda_{q+1}^{-1}$. If only a rough bound on $\ell_q$ is needed, then we will use $\ell_q\in(\lambda_q^{-2},\lambda_q^{-1}),$ which follows from the mild inequalities $\alpha<1/10$, $(b-1)\beta<1/5$ and $a^{3/5}>2.$   See also the definition of $\tau_q$ in \eqref{e:tau_q-and-t_i}.

We will use the notation $X\lesssim Y$ to denote that $X\le CY$ for some constant $C$ that depends on the various parameters $a,b,\alpha,\beta,\vinreg$, but not on $q$. By routine calculations similar to the above, we will freely use inequalities like $\lambda_{q+1}\lesssim \lambda_q^b$ and $\lambda_q^b \lesssim \lambda_{q+1}$, so long as we can lose some small positive power of $\lambda_q$ to absorb the implicit constants by taking $a\gg1$. The implicit constant may also depend on a certain number of derivatives $N$, but this is chosen to be a fixed number depending on $b$ and $\beta$ in Section \ref{s:est-stress}.
\subsection{Euler--Reynolds flow with smoothed initial data}
As is usual in convex integration schemes, we consider a modification of \eqref{e:euler} with a stress tensor error term $\RR_q$ that tends to $0$ in the sense of distributions. However,  we propagate an initial condition as well. That is, we require $(v_q,p_q,\RR_q)$ to solve
\begin{equation}
\left\{ \begin{alignedat}{-1}
\del_t v_q+\div (v_q\otimes v_q)  +\nabla p_q   &=  \div \RR_q,
 \\
  \nabla \cdot v_q &= 0,
  \\ v_q \big|_{t=0}&=\vin*\psi_{\ell_{q-1}}.
\end{alignedat}\right.  \label{e:subsol-euler}
\end{equation}
where
$a \otimes b\coloneq ab^\TT  = (a_i b_j)_{i,j=1}^3$, and the divergence of a 2-tensor $M=(M_{ij})_{i,j=1}^3$ is defined to be the vector $\div M$ with components (using the Einstein summation convention)
\begin{align*}
(\div M)_i \coloneq   \partial_j M_{ij}.
\end{align*}
In particular, since $v_q$ is divergence-free, $\div (v_q\otimes v_q) = (v_q \cdot\nabla )v_q$.
 The tensor $\RR_q$ is required to be a symmetric, trace-free $3\times3$ matrix, i.e.
\begin{gather}
     \RR_q = \RR_q^\TT ,  \qquad \tr\RR_q = \sum_{i=1}^3 (\RR_q)_{ii} =  0. \label{e:RR-cond}
\end{gather}
We enforce mean zero conditions on $v_q$ and $p_q$:
 \begin{align}
     \int_{\mathbb T^3} v_q  \dd x= 0, \qquad \int_{\mathbb T^3} p_q \dd x = 0. \label{e:mean-zero}
 \end{align}
\subsection{The main iterative proposition}
\label{ss:main-iterative-prop}
We will choose  $a\gg1$, and $b>1$ very close to 1.
We then choose a smooth function $e:[0,\infty]\to (\|\vin\|_{L^2}^2,\infty)$ such that
\begin{align}
t\in[1-\tau_0,T]\implies \delta_{2}\lambda_1^{-\alpha} &\le e(t) -\int_{\mathbb T^3} |\vin|^2 \dd x \le \delta_{2} .
    \label{e:vin-energy-estimate}
\end{align}

The estimates we propagate inductively are:
\begin{align}
    \|v_q\|_{0} &\le 1-\delta_q ,
    \label{e:vq-C0}
    \\
    \|v_q\|_{1} &\le M \delta_{q+1}^{1/2} \lambda_q ,
    \label{e:vq-C1}
    \\
    \|\RR_q\|_{0} &\le \delta_{q+1}\lambda_q^{-3\alpha} ,
    \label{e:RR_q-C0}
    \\
   t\in[ 1-\tau_{q-1},T] \implies \delta_{q+1}\lambda_q^{-\alpha} &\le e(t) -\int_{\mathbb T^3} |v_q(x,t)|^2 \dd x \le \delta_{q+1} .
    \label{e:energy-q-estimate}
\end{align}
The constant $M$ in \eqref{e:vq-C1} is universal (in particular independent of $q$): it is defined in \eqref{d:M}.

Note that the above estimates cannot be satisfied when using arbitrary initial data $\vin$, but this restriction is easily removed in the proof of Theorem \ref{t:main} by a scaling argument.
\begin{prop}
\label{p:main-prop} Let $M$ be the universal constant defined in \eqref{d:M}.  Let $T\in[1,\infty]$, $0<\beta<\vinreg <1/3$, and let $\vin\in C^{\vinreg}(\mathbb T^3)$ be a  periodic function with zero mean,  such that $\|\vin\|_{C^{\vinreg}}\le 1$, and let $e:[0,T]\to[0,\infty)$ be any smooth function.
Let $b\in(1,b_0)$ where $b_0=b_0(\beta,\vinreg)$ is given by the inequalities  \eqref{e:params1-b-alpha}, \eqref{e:params3-b}, and \eqref{e:params5-b-alpha},
\begin{align*}
    b_0 \coloneq 1+\min\left(1, \frac{\vinreg - \beta}{4\beta+3}, \frac{1-3\beta}{2\beta} \right).
\end{align*}
There exist  $\alpha_0$ and $a_0$ such that for $\alpha<\alpha_0(\beta,\vinreg,b)$ and $a>a_0(\beta,\vinreg, b,\alpha)$, the following holds.
Let $(v_q,p_q,\RR_q)$ solve
\eqref{e:subsol-euler} and satisfy \eqref{e:vq-C0}--\eqref{e:energy-q-estimate}.
Then there exists smooth functions $(v_{q+1}, p_{q+1}, \RR_{q+1})$, satisfying \eqref{e:subsol-euler},
\eqref{e:vq-C0}--\eqref{e:energy-q-estimate}
with $q$ replaced by $q+1$, and such that
\begin{align}
        \|v_{q+1} - v_q\|_{0} +\frac1{\lambda_{q+1}} \|v_{q+1}-v_q\|_1 &\le  M \delta_{q+1}^{1/2}.
        \label{e:velocity-diff}
\end{align}
The upper bound $\alpha_0$ is given by the inequalities \eqref{e:params0}, \eqref{e:params1-b-alpha}, \eqref{e:params2-alpha}, \eqref{e:params3-b}, \eqref{e:params6-alpha}, and  \eqref{e:params5-b-alpha},
 and $a_0$ is taken sufficiently large to remove various universal constants throughout the argument. 
\end{prop}

\subsection{Proof of Theorems \ref{t:main} and \ref{t:quant} using Proposition \ref{p:main-prop}}
\label{ss:proof-of-main}
Let $0<\beta<\vinreg<1/3$ and choose $b>b_0(\beta,\vinreg)$, $\alpha<\alpha_0(\beta,\vinreg,b)$ and $a>a_0(\beta,\vinreg,b,\alpha)$ as in Proposition \ref{p:main-prop}. We will later impose some further mild conditions on $a$ and $b$.

To start the iteration, we define $(v_1,p_1,\RR_1)$ by	
\begin{gather*}
    v_1\coloneq \vin*\psi_{\ell_{0}},
    \qquad  p_1\coloneq |v_1|^2-\int_{\mathbb T^3} |v_1|^2\dd x,
    \\ \RR_1\coloneq v_1\otimes v_1 -|v_1|^2 I_{3\times 3}.
\end{gather*}
 It is easy to check that they solve \eqref{e:subsol-euler}, \eqref{e:RR-cond}, and \eqref{e:mean-zero}.  In addition, we first assume that  they satisfy $\|\vin\|_{C^{\vinreg}}\le 1$, \eqref{e:vq-C0}--\eqref{e:RR_q-C0}.
We take a smooth function $e:[0,T]\to[0,\infty)$ satisfying
\begin{align*}
      t\in[1-\tau_0,T]\implies e(t) \in (\|{\vin}\|_{L^2}^2+\delta_2\lambda_1^{-\alpha} , \|{\vin}\|_{L^2}^2 + \delta_2 )
\end{align*} which therefore satisfies \eqref{e:vin-energy-estimate}. As $[1-\tau_0,T]$ is compact, there exists a constant $c=c(\vin)\in(1,2)$ such that
 \begin{align*}
      t\in[1-\tau_0,T]\implies e(t)\in (\|{\vin}\|_{L^2}^2+c\delta_2\lambda_1^{-\alpha} , \|{\vin}\|_{L^2}^2 + \delta_2/c ) .
\end{align*}
 In order to ensure that \eqref{e:energy-q-estimate} holds for $q=1$, we need to ensure that
$ \left| \int_{\mathbb T^3} |\vin|^2 \dd x - \int_{\mathbb T^3} |v_1 |^2 \dd x \right| \ll \delta_2.$
For this, we use the  Constantin--E--Titi estimate \cite{zbMATH00722097}:
\begin{align*}
    \left| \int_{\mathbb T^3} |\vin|^2 \dd x - \int_{\mathbb T^3} |v_1 |^2 \dd x \right| &\le \left| \int_{\mathbb T^3} |\vin|^2*\psi_{\ell_{0}} \dd x - \int_{\mathbb T^3} |\vin *\psi_{\ell_{0}} |^2 \dd x \right|
    \\ & \lesssim \|  |\vin|^2*\psi_{\ell_{0}}-  |\vin *\psi_{\ell_{0}} |^2 \|_{C^0}
    \\& \lesssim \ell_{0}^{2\vinreg}\|\vin\|_{C^{\vinreg}}^2 \ll \lambda_0^{-2\vinreg},
\end{align*}
which can be made smaller than $\delta_2\sim \lambda_0^{-2\beta b^2}=\lambda_0^{-2\beta (b-1)^2-4 \beta (b-1) -2\beta}$ by further restricting $b-1<\frac{\vinreg-\beta}{4\beta}$ and then increasing $a$ (at the beginning of the proof) to absorb the implicit constants, if necessary. Then, Proposition \ref{p:main-prop} applies inductively, giving rise to a $C^0$ convergent sequence of functions $v_q\to v$,  that solve \eqref{e:euler}, with $\| v\|_{L^2}^2(t)= e(t)>\|\vin\|_{L^2}^2$ for $t\in[1,T]$. A standard interpolation argument with the $C^1$ estimate in \eqref{e:velocity-diff} shows that $v\in C^{\beta}_{x,t}$ (see \cite{zbMATH07038033}).

To allow initial data of arbitrary size, we apply the theorem to
\[\widetilde \vin \coloneq  \frac{\vin}{\Gamma}, \qquad \Gamma \coloneq \max\left( \frac{\|\vin\|_{C^{\vinreg} }}{1-\delta_1}, \frac{\|\vin\|_{C^0}}{M\delta_2^{1/2}\lambda_1^{-3\alpha/2}} \right).  \]
The scaling is chosen precisely so that \eqref{e:vq-C0} and \eqref{e:RR_q-C0} hold. As before,  we set $v_1=\widetilde{\vin}*\psi_{\ell_0}$. Note that (since $\lambda_0<\lambda_1$ and $1+3\alpha/2-\beta>0$)
\[ \|v_1\|_1\le \frac1\Gamma \|\vin\|_{C^0}\ell_0^{-1} \le M \ell_0^{-1}  \delta_2^{1/2}\lambda_1^{-3\alpha/2} \le M \delta_2^{1/2} \lambda_1, \]
 and hence \eqref{e:vq-C1} holds.  We then take an arbitrary choice of $\tilde e:[0,\Gamma T]\to[0,\infty)$ that satisfies \eqref{e:vin-energy-estimate} with $\Gamma  T$ in place of $T$, and construct a solution $\widetilde v$ on the time interval $[0,\Gamma T]$, with $\|\widetilde v\|_{L^2}^2(t)= \tilde e(t)>\|\widetilde{\vin}\|_{L^2}^2$ for $t\in[1,\Gamma T]$.  Then we define
\[ v(x,t) \coloneq \Gamma\widetilde v(x,\Gamma t), \quad \text{ and }\quad e(t) \coloneq \Gamma^2 \tilde e(\Gamma t).\]
Due to the symmetry of the Euler equations, $v$ also solves Euler, but on the time interval $[0,T]$, and with the initial data $\vin$. In addition, since $1>\Gamma^{-1}$, we also obtain that $\|v\|_{L^2}^2(t)=e(t)>\|\vin\|_{L^2}$ for times $t\ge 1$.
\hfill $\qedsymbol$\\[0.5em]

A closer inspection of the proof gives Theorem \ref{t:quant}: for an initial data $\vin\in C^{\vinreg}$, we can choose for times $t\ge1$ an energy profile with image in $(\|\vin\|_{L^2}^2 + \Gamma^2\delta_2\lambda_1^{-\alpha}, \|\vin\|_{L^2}^2 + \Gamma^2\delta_2)$. Since for a fixed $\vin$, $\|\vin\|_{C^{\vinreg}} / (1-\delta_1)$ is bounded uniformly in the parameters $a,b,\alpha$, we should generically have  when $a\gg1$ that $\Gamma^2=\frac{\|\vin\|_0^2}{M\delta_2\lambda_1^{-3\alpha}}$. In such a situation, we would have
    \[ M\|\vin\|_{C^0}^2\lambda_1^{2\alpha} \le e(t) - \int_{\mathbb T^3}|\vin|^2 \dd x \le M\|\vin\|_{C^0}^2\lambda_1^{3\alpha}.\]
In particular, there is a large lower bound (as $\alpha>0$ and $a\gg 1$), and in order to achieve larger energies, we need the lower bound to increase as well.

The remainder of the paper is devoted to the proof of Proposition \ref{p:main-prop}.
\subsection{Proof sketch for Proposition \ref{p:main-prop}}
\label{ss:sketch-pf-main-prop}

Starting from a tuple $(v_q,p_q,\RR_q)$ satisfying the estimates as in Proposition \ref{p:main-prop}, the broad scheme of the iteration is as follows.
\begin{enumerate}
    \item $(v_{\ell_q},p_{\ell_q},\RR_{\ell_q})$ are defined by mollification. This step differs from previous works because we include a mollified initial data in the iteration scheme.
    \item Next we define a family of exact solutions to Euler $(\vex_i,\pex_i)$, $i\ge1$, by exactly solving the Euler equations with smooth initial data $v\big|_{t=t_i}= v_{\ell_q}(t_i)$, where $t_i=i\tau_q$ defines an evenly spaced  paritition of $[0,T]$.  For  $i=0$, we improve the initial data from $v_\ell \big|_{t=0}=\vin *\psi_{\ell_{q-1}}*\psi_{\ell_{q}}$ to $\vex_0\big|_{t=0}=\vin * \psi_{\ell_{q}}$, but this creates a small mismatch. For the estimates to work at this step, we require $\vinreg>\beta$. \label{step}
    \item These solutions are glued together using a partition of unity, resulting in the tuple $(\vv_q,\pp_q,\RRR_q)$. The stress error term is zero when $\vv_q=\vex_i$ for some $i\ge0$.
    \item We define $v_{q+1}$ by constructing a perturbation $w_{q+1}$ and setting $v_{q+1} = \vv_q+w_{q+1}$. In order to achieve the optimal regularity, we construct some key cutoffs functions $\eta_i$ (see \eqref{e:defn-cutoff-total}) and use the Mikado flows from \cite{zbMATH06710292} to create the principal part $\wpq$ of the perturbation. We also define an incompressibility corrector term $\wcq$ so that  $w_{q+1} \coloneq \wpq+\wcq$ is divergence-free.
    \item Once we have fixed the definition of $v_{q+1}$, the term $\div \RR_{q+1}$ is fixed since $(v_{q+1},\RR_{q+1})$ has to solve \eqref{e:euler}. This leads to the natural definition of $\RR_{q+1}$ using the inverse divergence operator $\mathcal R$.
    \item We prove that the inductive estimates \eqref{e:vq-C0}--\eqref{e:energy-q-estimate} hold with $q$ replaced with $q+1$.
    \end{enumerate}
Step 4 is moderately involved, and breaks into the following sub-steps, where we have combined the approach of \cite{zbMATH06976275} and \cite{zbMATH07038033} in order to achieve control of the energy:
\begin{enumerate}
    \item[4a.] For times $t\ge 1$, we use the `squiggling' cutoffs $\eta_i$ from \cite{zbMATH07038033} that allow energy to be added at such times, even outside the support of $\RRR_q$, while cancelling a large part of $\RRR_q$ norm.
    \item[4b.] For times $t< 1$, we instead use the straight cutoffs of \cite{zbMATH06976275}, so that we can ensure we do not adjust the solution near $t=0$. We leave some space in order to transition between the two approaches. Together with the modification of Step \ref{step}, this ensures that  we  have $v_{q+1}\big|_{t=0} = \vv_q\big|_{t=0} + w_{q+1}\big|_{t=0} = \vin*\psi_{\ell_q}$, so that the initial data is propagated inductively.
     \item[4c.] Then we create the principal part $\wpq$ of the perturbation using Mikado flows and the back-to-labels map of $\vv_q$, and add a small corrector term $\wcq$  to enforce the incompressibility constraint.
\end{enumerate}


\section{The iterative step}
In this section, we give some details for the construction of $(v_{q+1},p_{q+1},\RR_{q+1})$ from $(v_{q},p_{q},\RR_{q})$. Our construction is a modification of \cite{zbMATH07038033}, so we direct the reader there at a number of points for  details that are not elaborated here.
\label{s:construct}

\subsection{Mollification}
\label{ss:mollification}

We recall that $\ell_q$ is defined in \eqref{e:ell}. The functions $v_{\ell_q}$ and $\RR_{\ell_q}$ are defined with the spatial mollifier \eqref{e:defn-mollifier-x},
\begin{align*}
    v_{\ell_q} &\coloneq v_q * \psi_{\ell_q},
 &   \RR_{\ell_q} &\coloneq \RR_q * \psi_{\ell_q}  - (v_q \ootimes v_q) * \psi_{\ell_q}  + v_{\ell_q} \ootimes v_{\ell_q} , \label{e:v_ell}
\end{align*}
where for vectors $a,b\in\mathbb R^3$, we have written  $a\ootimes b \coloneq a\otimes b - (a\cdot b)I_{3\times 3}$, which is a trace-free version of $a\otimes b$.  $v_{\ell_q}$ is also divergence-free, and $\RR_{\ell_q}$ is also symmetric and trace-free. They solve the initial value problem
\begin{equation*}
\left\{ \begin{alignedat}{-1}
\del_t v_{\ell_q} +\div (v_{\ell_q}\otimes v_{\ell_q})  +\nabla p_{\ell_q}   &=  \div \RR_{\ell_q} ,
\\
  \nabla \cdot v_{\ell_q} &= 0,
  \\
  v_{\ell_q}\big|_{t=0} &= \vin*\psi_{\ell_{q-1}}*\psi_{\ell_{q}},
\end{alignedat}\right.  \label{e:mollified-euler}
\end{equation*}
where $p_{\ell_q} \coloneq p_q *\psi_{\ell_q} -|v_q|^2 + |v_{\ell_q}|^2,$ and we have used the identity $\div(fI_{3\times 3}) =\nabla f$ for scalar fields $f$. Standard mollification estimates and the quadratic Constantin--E--Titi estimate \cite{zbMATH00722097} give the following proposition.
\begin{prop}[Estimates for mollified functions]\label{p:estimates-for-mollified}
\begin{align}
\|v_{\ell_q}-v_{q}\|_{0} &\lesssim \delta_{q+1}^{1 / 2}  \lambda_{q}^{-\alpha} \label{e:v_ell-vq}
\\
\|v_{\ell_q}\|_{N+1} &\lesssim \delta_{q}^{1 / 2} \lambda_{q} \ell_q^{-N} && \forall N \geq 0, \label{e:v_ell-CN+1}
\\
\|\RR_{\ell_q}\|_{N+\alpha} &\lesssim \delta_{q+1} \ell_q^{-N+\alpha} && \forall N \geq 0 ,\label{e:R_ell}
\\
\Big|\int_{\mathbb{T}^{3}} | v_{q}|^{2}-|v_{\ell_q}|^{2} \dd x \Big| &\lesssim \delta_{q+1} \ell_q^{\alpha}  && \forall t\in[0,T].\label{e:energy-v_ell}
\end{align}
\end{prop}
\begin{proof}
See \cite[Proposition 2.2]{zbMATH07038033}.
\end{proof}

\subsection{Classical Exact  flows}
\label{ss:exact}
We define  $\tau_q\in(0,\ell_q^{\alpha})$ and the  sequence of initial times $t_i$ ($i\in\mathbb Z_{\ge0}$) by
 \begin{align}
    \tau_q \coloneq \frac{\ell_q^{2\alpha}}{\delta_q^{1/2} \lambda_q}, \qquad t_i\coloneq i\tau_q.   \label{e:tau_q-and-t_i}
\end{align}
Note $\|v_{\ell_q}\|_{1}\tau_q \lesssim \ell_q^{2\alpha}$, so $\|v_{\ell_q}   \|_{1} \le \tau_q^{-1}.$ That is, $\tau_q$ is chosen such that
\begin{align*}
    \delta_{q+1}^{1/2}\tau_q \ell_q^{-1} = \ell_q^{2\alpha}\lambda_q^{3\alpha/2} \le \lambda_q^{-\alpha/2} \le 1. 
\end{align*}
We define  $(\vex_i,\pex_i)$ for $i\ge 1$ to be the unique (exact) solutions
to the Euler equations, defined on some interval around $t_i$, with initial data
 $\vex_i\big|_{t=t_i} = v_{\ell_q}(t_i)$. That is,
  \begin{align*}
 \left\{ \begin{alignedat}{-1}
\del_t \vex_i +(\vex_i\cdot\nabla) \vex_i  +\nabla \pex_i   &=  0,
\\
  \nabla \cdot \vex_i  &= 0,
  \\
  \vex_i\big|_{t=t_i} &= v_{\ell_q}(\cdot,t_i)
\end{alignedat}\right. & 
 \intertext{
For $i=0$, we define $(\vex_0,\pex_0)$ to be the classical solution to the Euler equations starting from $\vin*\psi_{\ell_q}$:
}
\left\{ \begin{alignedat}{-1}
\del_t \vex_0 +(\vex_0\cdot\nabla) \vex_0  +\nabla \pex_0   &=  0,
\\
  \nabla \cdot \vex_0  &= 0,\phantom{\!(\cdot,t_i)}
  \\
  \vex_0\big|_{t=0} &= \vin * \psi_{\ell_q}
\end{alignedat}\right.  &
 \end{align*}
 We note that there is a mismatch of initial data: $\vex_0\big|_{t=0}\neq v_{\ell_q}\big|_{t=0} = \vin*\psi_{\ell_q}*\psi_{\ell_{q-1}}.$
 This difference needs to be estimated below, and relies on the parameter inequality $\beta<\vinreg$.
\begin{prop}[Estimates for classical exact solutions to Euler {\cite[Prop 3.1]{zbMATH07038033}}]
\label{p:exact-euler}
Let $\tilde N\ge1$ and $\tilde\alpha\in(0,1)$. The unique $C^{\tilde N+\tilde\alpha}$ solution $V$ to \eqref{e:euler} with initial data $V_0\in C^{\tilde N+\tilde\alpha}$ is defined at least for $t\in [-\tilde T,\tilde T]$, where $\tilde T= \frac c{\|V_0\|_{{1+\tilde\alpha}}}$ for some universal $c>0$, and satisfies the following derivative estimates for $0\le N\le \tilde N$,
\[ \|V\|_{C^0([-\tilde T,\tilde T];C^{N+\alpha})} \lesssim  \| V_0 \|_{{N+\alpha}} . \]

\end{prop}
Note that the definition of $\tau_q$ was chosen so that $\vex_i$ is well-defined on $[t_{i-1},t_{i+1}]$, once $a\gg1$ is sufficiently large.
\begin{prop}[Stability]
\label{p:stability}
Let $\matd v\coloneq \del_t+v\cdot\nabla$ denote the material derivative. Suppose $\beta,\alpha$ satisfy the constraints \eqref{e:params1-b-alpha}, \eqref{e:params2-alpha} below. For $i\ge 0$, $|t-t_i|\lesssim  \tau_q$, and $N\ge 0$, we have
\begin{align}
    \|\vex_i -v_{\ell_q} \|_{{N+\alpha}} &\lesssim \tau_q \delta_{q+1} \ell_q^{-N-1+\alpha} \label{e:stability-v},\\
    \|\nabla \pex_i -\nabla p_{\ell_q} \|_{{N+\alpha}} &\lesssim \delta_{q+1} \ell_q^{-N-1+\alpha},\label{e:stability-p}\\
    \|\matd{v_{\ell_q}}(\vex_i - v_{\ell_q}) \|_{{N+\alpha}} &\lesssim \delta_{q+1} \ell_q^{-N-1+\alpha}\label{e:stability-matd}.
\end{align}
\end{prop}
\begin{proof}
 The proof is similar to that of  \cite[Prop. 3.3]{zbMATH07038033}, but we need to treat $\vex_0$ separately.
 The equation for the difference is
 \begin{align} \del_t (\vex_i - v_{\ell_q}) + v_{\ell_q} \cdot\nabla (\vex_i - v_{\ell_q}) + (\vex_i-v_{\ell_q})\cdot\nabla \vex_i + \nabla (\pex_i - p_{\ell_q}) = - \div \RR_{\ell_q}.  \label{e:eqn-for-vex_i-v_ell}\end{align}
Taking the divergence of the above equation, we obtain a Poisson equation for $\pex_i-p_{\ell_q}$, leading to the estimates ($N\ge 0$, $i\ge1$)
\begin{align}
     & \|\nabla \pex_i{-}\nabla p_{\ell_q}\|_{{N+\alpha}}\notag  \\
     &= \|\Delta^{-1}\nabla \div \big[ v_{\ell_q} {\cdot}  \nabla (\vex_i{-}v_{\ell_q})+(\vex_i{-}v_{\ell_q}){\cdot} \nabla \vex_i +\div\RR_{\ell_q}\big] \|_{{N+\alpha}}  \notag
     \\
     &\lesssim   \|\nabla \vex_i\|_{{N+\alpha}} \|\vex_i-v_{\ell_q}\|_{{\alpha}} + \|\nabla \vex_i \|_{{\alpha}} \|\vex_i-v_{\ell_q}\|_{{N+\alpha}} + \| \RR_{\ell_q} \|_{{N+1+\alpha}} \notag \\
     &\quad +\|\nabla v_{\ell_q}\|_{{N+\alpha}} \|\vex_i-v_{\ell_q}\|_{{\alpha}} + \|\nabla v_{\ell_q}\|_{{\alpha}} \|\vex_i-v_{\ell_q}\|_{{N+\alpha}}  \notag \\
     &\lesssim \delta_q^{1/2}\lambda_q\ell_q^{-N-\alpha}\|\vex_i-v_{\ell_q}\|_{{\alpha}} + \delta_q^{1/2}\lambda_q\ell_q^{-\alpha}\|\vex_i-v_{\ell_q}\|_{{N+\alpha}} +\delta_{q+1}\ell_q^{-N-1+\alpha}, \label{e:nabla-p}
\end{align}
where we have used \eqref{e:v_ell-CN+1}, \eqref{e:R_ell}, Schauder estimates \eqref{e:schauder}, the identity $\nabla \div( (V \cdot\nabla) W )= \nabla \div( (W \cdot\nabla) V)$ to move the gradient off of $\vex_i-v_{\ell_q}$, and the estimate   $\|\vex_i\|_{{N+\alpha}} \le\|v_{\ell_q}\|_{{N+\alpha}}$ from Proposition \ref{p:exact-euler}.

For $i=0$, Proposition \ref{p:exact-euler} instead gives  $\|\vex_0\|_{{N+\alpha}} \le\|\vin*\psi_{\ell_q}\|_{{N+\alpha}}$.
Since $\|\vin\|_{\vinreg} \le 1$, we get by \eqref{e:convolution-holder-estimate} that
\[ \|\vin*\psi_{\ell_q} \|_{N+1+\alpha} \lesssim \|\vin\|_{\vinreg} \ell_q^{-N-1-\alpha+\vinreg} \lesssim \ell_q^{-N-1-\alpha+\vinreg}. \]
In order to  have \eqref{e:nabla-p} for $i=0$, we therefore require
\begin{align}
    \ell_q^{-N-1-\alpha+\vinreg} \lesssim \delta_q^{1/2} \lambda_q \ell_q^{-N-\alpha}
     \iff \ell_q^{\vinreg-\beta }\lesssim (\lambda_q\ell_q)^{1-\beta}. \label{e:required-ineq1}
\end{align}
For this, as $\lambda_q\ell_q<1$, it suffices to make $\lambda_q^{-(\vinreg-\beta)}\lesssim \lambda_q\ell_q=\delta_{q+1}^{1/2}\delta_q^{-1/2}\lambda_q^{-3\alpha}$, and this is true if we impose 
\begin{align}
    b-1\le \min\left(1,\frac{\vinreg - \beta}{4\beta}\right),\quad\text{and}\quad  3\alpha<\beta(b-1)     \label{e:params1-b-alpha},
\end{align}
since this implies
\[ \lambda_q^{-(\vinreg-\beta)} \le \lambda_q^{-4\beta (b-1)} \ll \lambda_q^{-\beta(b-1) - 3\alpha } \le 2 \delta_{q+1}^{1/2}\delta_q^{-1/2}\lambda_q^{-3\alpha}. \] 
Hence, from \eqref{e:eqn-for-vex_i-v_ell} we obtain for $N=0$, using Proposition \ref{p:estimates-for-mollified}
\begin{align*} \|\matd{v_{\ell_q}} (\vex_i - v_{\ell_q}) \|_{\alpha}  &\lesssim   \delta_q^{1/2}\lambda_q\ell_q^{-\alpha} \|\vex_i-v_{\ell_q}\|_{\alpha} + \|\RR_{\ell_q} \|_{{1+\alpha}}\\ &  \lesssim   \tau_q^{-1} \|\vex_i-v_{\ell_q}\|_{\alpha} + \ell_q^{-1+\alpha} \delta_{q+1}.
\end{align*}
From standard estimates \cite[Appendix D]{zbMATH06456007} for the transport equation, it follows that
\begin{align*}
&\|v_{\ell_q} - \vex_i\|_{\alpha} (t)
\lesssim \|v_{\ell_q} - \vex_i\|_{\alpha}(t_i) + \tau_q \delta_{q+1} \ell_q^{-1+\alpha}+\int_{t_i}^t  \tau_q^{-1}\|v_{\ell_q} -\vex_i\|_{\alpha}(s)  \dd s .
\end{align*}
For $i\ge1$, $v_{\ell_q}\big(t_i) = \vex_i\big(t_i)$. For $i=0$, it suffices to control the initial data by $\tau_q \delta_{q+1} \ell_q^{-1+\alpha}$. By \eqref{e:convolution-diff-holder-estimate} and $\|\vin\|_\beta \le 1$, we have
\begin{align*}
\|v_{\ell_q} - \vex_0\|_{\alpha}\big|_{t=0} = \|\vin*\psi_{\ell_q}-\vin*\psi_{\ell_q}*\psi_{\ell_{q-1}}\|_\alpha &\lesssim \|\vin\|_{\vinreg}\ell_{q-1}^{\vinreg-\alpha},
\end{align*}
and we need $\ell_{q-1}^{\vinreg-\alpha}\le\tau_q\delta_{q+1}\ell_q^{-1+\alpha}$, i.e. $ \ell_{q-1}^{\vinreg} \le \lambda_{q+1}^{-\beta}\ell_q^{3\alpha}\lambda_q^{3\alpha/2}\ell_{q-1}^\alpha. $
It suffices to take $b-1\le \min(1,\frac{\vinreg - \beta}{4\beta})$ as in \eqref{e:params1-b-alpha}, and
\begin{align}
 8 \alpha < \frac{\beta(b-1)}{b} \label{e:params2-alpha},
\end{align}
so that $-\vinreg\le - 4\beta(b-1)-\beta$, and hence (since $b<2$ and $\lambda_q^{-2}\le\ell_q$)
\[\lambda_{q-1}^{-\vinreg}  \lesssim  \lambda_{q-1}^{-\beta b^2}\lambda_q^{-\beta(b-1)/b} \le \lambda_{q-1}^{-\beta b^2} \lambda_q^{-8\alpha} \le \lambda_{q-1}^{-\beta b^2} \ell_q^{4\alpha}  \ll   \lambda_{q+1}^{-\beta} \lambda_q^{3\alpha/2}\ell_q^{3\alpha} \ell_{q-1}^{\alpha}, \]
which implies the wanted estimate.
Hence by \Gronwall{}'s inequality,
\[\| \vex_i-v_{\ell_q}\|_{\alpha} \le \tau_q\delta_{q+1}\ell_q^{-1+\alpha} \ee^{\tau_q^{-1}(t-t_i)}\lesssim \tau_q \delta_{q+1}\ell_q^{-1+\alpha}, \]
which is \eqref{e:stability-v} for $N=0$ and all $i\ge0$. The estimates for $N>0$ are similar:
Let $\sigma$ be a multiindex with $|\sigma|=N$. Then we write\[ \matd{v_{\ell_q}} D^\sigma (\vex_i-v_{\ell_q}) = [\matd{v_{\ell_q}},D^\sigma](\vex_i-v_{\ell_q}) + D^\sigma \matd{v_{\ell_q}} (\vex_i-v_{\ell_q}). \]
The commutator term cancels the highest derivative on $\vex_i-v_{\ell_q}$:
\begin{align}
    \| [\matd{v_{\ell_q}},D^\sigma](\vex_i-v_{\ell_q}) \|_{\alpha}
    \lesssim \|v_{\ell_q}\|_{N+\alpha}\|\vex_i-v_{\ell_q}\|_{1+\alpha} + \|v_{\ell_q}\|_{1+\alpha}\|\vex_i-v_{\ell_q}\|_{N+\alpha}\notag
    \\
\lesssim \|v_{\ell_q}\|_{N+1+ \alpha}\|\vex_i-v_{\ell_q}\|_{\alpha} + \|v_{\ell_q}\|_{1+\alpha}\|\vex_i-v_{\ell_q}\|_{N+\alpha}\notag
    \\
    \lesssim \delta_{q}^{1 / 2} \lambda_{q} \ell_q^{-N-\alpha}
\cdot  \tau_q \delta_{q+1} \ell_q^{-1+\alpha} + \delta_{q}^{1 / 2} \lambda_{q} \ell_q^{-\alpha}\|\vex_i-v_{\ell_q}\|_{N+\alpha}. \notag  \end{align}
(We used \Holder{} interpolation in the second inequality.) For the other term $D^\sigma \matd{v_{\ell_q}} (\vex_i-v_{\ell_q})$, we use the equation \eqref{e:eqn-for-vex_i-v_ell} to write
  \begin{align*} \|D^\sigma \matd{v_{\ell_q}} (\vex_i-v_{\ell_q})\|_\alpha  = &\lVert D^\sigma (\div \RR_{\ell_q} + \nabla  (\pex_i-p_{\ell_q}) - (\vex_i-v_{\ell_q})\cdot\nabla \vex_i)
    \|_{\alpha}.
\end{align*}
 With the estimates  \eqref{e:R_ell},  \eqref{e:nabla-p}, \eqref{e:required-ineq1} and  the estimate
\begin{align}
    \|D^\sigma ((\vex_i-v_{\ell_q})\cdot\nabla \vex_i)\|_\alpha
    \le \| \vex_i-v_{\ell_q}\|_{\alpha}\|\vex_i\|_{N+1+\alpha}+\| \vex_i-v_{\ell_q}\|_{N+\alpha}\|\vex_i\|_{1+\alpha}  \notag
    \\
    \le \tau_q \delta_{q+1} \ell_q^{-1+\alpha} \cdot \delta_{q}^{1 / 2} \lambda_{q} \ell_q^{-N-\alpha} + \delta_{q}^{1 / 2} \lambda_{q} \ell_q^{-\alpha} \|\vex_i-v_{\ell_q}\|_{N+\alpha}, \notag
\end{align}
 we obtain the bound
 \begin{align*}
    \| \matd{v_{\ell_q}} D^\sigma (\vex_i-v_{\ell_q}) \|_\alpha &\lesssim  \tau_q \delta_{q+1}  \delta_{q}^{1 / 2} \lambda_{q} \ell_q^{-1-N} + \delta_{q}^{1 / 2} \lambda_{q} \ell_q^{-\alpha} \|\vex_i-v_{\ell_q}\|_{N+\alpha} \notag
    \\ &\quad   + \delta_{q+1} \ell_q^{-N-1+\alpha} \notag
    \\ &\le  \delta_{q+1} \ell_q^{-N-1+\alpha} + \tau_q^{-1} \|\vex_i-v_{\ell_q}\|_{N+\alpha}.
    \end{align*}
By  estimates for the transport equation again and summing over  $\sigma$, then applying \Gronwall{}'s inequality, we obtain
\begin{align*}
    \|\vex_i-v_{\ell_q}\|_{N+\alpha} \lesssim  \tau_q \delta_{q+1} \ell_q^{-N-1+\alpha} ,
\end{align*} which gives \eqref{e:stability-v} and \eqref{e:stability-matd} for $N> 0$, and hence   \eqref{e:stability-p} for all $N> 0$ and $i\ge0$.
\end{proof}

\subsection{Estimates for Vector potentials}
\label{ss:vector-potentials}
By the well-known Helmholtz decomposition for smooth functions on the torus, a divergence-free field $V$ of zero mean is in fact a curl, i.e. $V=\nabla\times Z$ for an incompressible field  $Z\eqcolon \mathcal BV$  called the vector potential of $V$. The operator $\mathcal B=(-\Delta)^{-1}\curl$ is the `Biot--Savart operator'. Here we note some estimates for the vector potential of $\vex_i$. These are used in estimating the stress error $\RRR_q$ that arises from gluing different Euler flows together: see Proposition \ref{p:estimate-RRRq}.
\begin{prop}
   Define\footnote{A minor difference from the iteration scheme in \cite{zbMATH07038033} is that we propogate the zero mean condition, and hence that $\mathcal B \vex_i$ is well-defined.} $z_i \coloneq \mathcal B \vex_i$. Let $\beta$ and $\alpha$ satisfy the constraints \eqref{e:params3-b} below. For $\left|t-t_{i}\right| \leq \tau_{q}$,
\begin{align}\|z_{i}-z_{i+1}\|_{N+\alpha} &\lesssim \tau_{q} \delta_{q+1} \ell_q^{-N+\alpha},\label{e:deriv-vector-potential}
\\
\|\matd{v_{\ell_q}} (z_{i}-z_{i+1})\|_{N+\alpha} &\lesssim \delta_{q+1} \ell_q^{-N+\alpha}. \label{e:matd-vector-potential} \end{align}
\end{prop}
\begin{proof}
For $N\ge1$, \eqref{e:deriv-vector-potential} directly follows from the analogous estimates \eqref{e:stability-v} for $\vex_i$, and the boundedness of the zero-order operator $\nabla\mathcal B$ on \Holder spaces. For $N=0$, we first prove the intermediate estimate \eqref{e:matd-vector-potential-aux} that will also lead to \eqref{e:matd-vector-potential}. One can check that \eqref{e:eqn-for-vex_i-v_ell} can be written in terms of $\tilde z_i=\mathcal B(\vex_i-v_{\ell_q})$ as\footnote{Here, $(a\times \nabla)^j b^k \coloneq \epsilon_{jkl}a^k\del_{x^l}b^k.$ The identities we need follow from $\nabla\cdot \tilde z_i=0$. They are
$v_{\ell_q}\cdot\nabla (\vex_i-v_{\ell_q}) = \curl( (v_{\ell_q}\cdot\nabla)\tilde z_i ) + \sum_{j,k=1}^3\del_{x^k}\Big((\tilde z_i\times \nabla )^j v_{\ell_q}^k\Big)$
and $((\vex_i-v_{\ell_q})\cdot\nabla)\vex_i= \sum_{j,k=1}^3\del_{x^k}\Big((\tilde z_i\times \nabla )^k \vex_i^j\Big)$. }
\[  \curl(\del_t \tilde z_i + v_{\ell_q} \cdot\nabla \tilde z_i ) = -\del_{x^k}\Big((\tilde z_i\times \nabla )^j v_{\ell_q}^k+(\tilde z_i\times \nabla )^k \vex_i^j\Big)-\nabla(\pex_i-p_{\ell_q})-\div\RR_{\ell_q}, \]
where there is an implicit sum over $j,k=1,2,3$.  On taking the curl and using the identity $\curl\curl=-\Delta+\nabla\div$, we find that $\matd{v_{\ell_q}} \tilde z_i$ solves the Poisson equation
\begin{align*}
     & {-\Delta}(\del_t \tilde z_i + (v_{\ell_q} \cdot\nabla) \tilde z_i ) \\ &= -\nabla\div((\tilde z_i\cdot\nabla)v_{\ell_q})- \curl\del_{x^k}\Big((\tilde z_i\times \nabla )^j v_{\ell_q}^k+(\tilde z_i\times \nabla )^k \vex_i^j\Big)-\curl \div \RR_{\ell_q} .
\end{align*}
 Schauder estimates for the Poisson equation and the boundedness of singular integrals on \Holder spaces therefore give
\begin{align}
    \| \matd{v_{\ell_q}} \tilde z_i\|_{N+\alpha} &\lesssim \|(\tilde z_i\cdot\nabla)v_{\ell_q}\|_{N+\alpha} \notag
    \\
    &\quad + \sum_{j,k=1}^3 \|(\tilde z_i\times \nabla )^j v_{\ell_q}^k+(\tilde z_i\times \nabla )^k \vex_i^j\|_{N+\alpha} + \|\RR_{\ell_q}\|_{N+\alpha} \notag
    \\& \lesssim (\|v_{\ell_q}\|_{N+1+\alpha} + \|\vex_i\|_{N+1+\alpha})\|\tilde z_i \|_{\alpha} \notag
    \\& \quad + (\|v_{\ell_q}\|_{1+\alpha} + \|\vex_i\|_{1+\alpha})\|\tilde z_i \|_{N+\alpha} + \delta_{q+1}\ell_q^{-N+\alpha}. \label{e:matd-vector-potential-aux}
\end{align}
In the case $i\neq0$, we can now easily prove both inequalities for all $N$. Then the $N=0$ case of \eqref{e:matd-vector-potential-aux} and standard estimates on the transport equation (since $\tilde z_i\big|_{t=0}=0$) give
\begin{align*}
    \|\tilde z_i\|_{\alpha} \lesssim \delta_q^{1/2}\lambda_q\ell_q^{-\alpha}\int_{t_i}^t\|\tilde z_i\|_{\alpha} \dd s + \delta_{q+1} \ell_q^{\alpha}|t-t_i|\lesssim \tau_q^{-1} \int_{t_i}^t\|\tilde z_i\|_{\alpha} \dd s + \delta_{q+1} \ell_q^{\alpha}\tau_q,\end{align*}
    and hence \Gronwall{}'s inequality proves \eqref{e:deriv-vector-potential} for $N=0$. Then \eqref{e:matd-vector-potential} follows from \eqref{e:matd-vector-potential-aux} for all $N\ge0$.

For $i=0=N$, the initial data for $\tilde z_0$ is
\[\tilde z_0\big|_{t=0} = \mathcal B( \vin * \psi_{\ell_q} - \vin * \psi_{\ell_{q-1}} * \psi_{\ell_q}) = (\mathcal B\vin)  * \psi_{\ell_q} -(\mathcal B\vin)*\psi_{\ell_{q-1}} * \psi_{\ell_q}. \]
For this we use \eqref{e:convolution-diff-holder-estimate2} to find
\[\|(\mathcal B\vin)  * \psi_{\ell_q} -(\mathcal B\vin)*\psi_{\ell_{q-1}} * \psi_{\ell_q}\|_{\alpha} \le [\nabla \mathcal B \vin * \psi_{\ell_q}]_{\vinreg} \ell_{q-1}^{1+\vinreg-\alpha} \lesssim \ell_{q-1}^{1+\vinreg-\alpha}. \]
To match with the estimate for $i>0$, we require
\[ \ell_{q-1}^{1+\vinreg - \alpha} \le \delta_{q+1}\ell_q^{\alpha} \tau_q \iff \ell_{q-1}^{\vinreg } \le  \delta_{q+1}^{1/2} \ell_q^{3\alpha} \lambda_q^{3\alpha/2} \frac{\ell_q}{\ell_{q-1}} \ell_{q-1}^\alpha \]
Since $\ell_{q+1}=\ell_{q}^b$ up to universal multiplicative constants, $\lambda_{q-1}^{-2(b-1)} \le \ell_{q-1}^{b-1} \lesssim \frac{\ell_q}{\ell_{q+1}}$ and it suffices to assume that
\begin{align}
     b-1<\min\left(1, \frac{\vinreg-\beta}{4\beta+3}\right)\text{ and } \quad  8\alpha \le \frac{b-1}{b}  \label{e:params3-b} ,
\end{align}
so that
\[ \lambda_{q-1}^{-\vinreg}\le \lambda_{q-1}^{-\beta b^2 - (b-1)\beta - 3(b-1)} \ll \delta_{q+1}^{1/2} \lambda_{q-1}^{-2(b-1)} \lambda_q^{-(b-1)/b} \le  \delta_{q+1}^{1/2}\frac{\ell_q}{\ell_{q-1}} \lambda_q^{-8\alpha}.\]
This is enough because
 $\lambda_q^{-8\alpha} \le \ell_q^{4\alpha} \ll \ell_q^{3\alpha}\ell_{q-1}^\alpha \lambda_q^{3\alpha/2}$, giving the required inequality for $a\gg1$. This means that we can bound as before,
\begin{align*}
    \|\tilde z_0\|_{\alpha} &\lesssim \|\mathcal B\vin -(\mathcal B\vin)*\psi_{\ell_q} \|_\alpha + \tau_q^{-1} \int_{0}^t\|\tilde z_0\|_{\alpha} \dd s + \delta_{q+1} \ell_q^{\alpha}\tau_q\\
    &\lesssim \tau_q^{-1} \int_{0}^t\|\tilde z_0\|_{\alpha} \dd s + \delta_{q+1} \ell_q^{\alpha}\tau_q.\end{align*}
Using \Gronwall{} again gives \eqref{e:deriv-vector-potential}, and \eqref{e:matd-vector-potential} again follows from \eqref{e:matd-vector-potential-aux} for all $N\ge0$.
   \end{proof}
\subsubsection{Gluing exact  flows}
\label{ss:gluing}
This subsection is taken from \cite[Section 4]{zbMATH07038033}, so we omit some routine calculations and proofs: despite our modifications, we have shown that the same estimates hold for $\vex_i$, $v_q$, and $v_{\ell_q}$ as in that paper, so the proofs for the results here need no modification.

Define the intervals $I_i,J_i$ ($i\ge -1)$ by $    I_i \coloneq [ t_i + \frac{\tau_q}3,\ t_i+\frac{2\tau _q}3]$ ,
and
$   J_i \coloneq (t_i - \frac{\tau_q}3,\ t_i+\frac{\tau _q}3).$
They partition $\mathbb R$. Define $i_{\textup{max}}$ to be the smallest number so that
$[0,T]\subseteq J_0 \cup I_0 \cup J_1 \cup I_1 \cup \dots \cup J_{i_{\textup{max}}} \cup I_{i_{\textup{max}}},$ i.e.
  \[i_{\textup{max}}\coloneq\sup\{ i\ge0 :   (J_i\cup I_i)\cap [0,T]\neq \emptyset \}\le \left\lceil \frac T{\tau_q}\right\rceil.\]
Also let $\{ \chi_i\}_{i=0}^{i_{\textup{max}}}$  be a partition of unity such that for all $N\ge 0$,
\begin{align}
     \supp \chi_i=I_{i-1} \cup J_i \cup I_{i}, \quad  \chi_i |_{J_i} =1, \quad  \|  \del_t^N \chi_i\|_{C^0_t} \lesssim \tau_q^{-N} . \label{e:chi_i-properties}
\end{align}
In particular, for $|i-j|\ge2$, $\supp \chi_i \cap \supp \chi_{j} =\emptyset$. By a slight abuse of notation, we write $ \chi_i \vex_i \coloneq 0$ for $t$ outside the support of $\chi_i$ even when $\vex_i$ is not defined, and similarly for other functions.  Then we define the glued velocity and pressure $(\vv_q,\ppp_q)$ by
\begin{align*}
    \vv_q(x,t) &\coloneq \sum_{i=0}^{i_{\textup{max}}} \chi_i(t) \vex_i(x,t) , 
    \\
    \ppp_q(x,t) &\coloneq \sum_{i=0}^{i_{\textup{max}}} \chi_i(t) \pex_i(x,t). 
\end{align*}
The definition \eqref{e:chi_i-properties} of $\chi_i$ implies that $\vv_q$ is still divergence-free, and \begin{align*}
    t\in I_i\implies\vv_q(x,t)&=\vex_i(x,t)\chi_i(t) + \vex_{i+1}(x,t)\chi_{i+1}(t),\\
    t\in J_i\implies \vv_q (x,t)&=\vex_i(x,t).
\end{align*}
Hence, $\vv_q$ is an exact Euler flow for the times $t\in J_i$.
For $t\in I_i$, we define
\begin{align*}
    \RRR_q &\coloneq
        \del_t \chi_i \mathcal R(\vex_i-\vex_{i+1} ) - \chi_i(1-\chi_i)(\vex_i-\vex_{i+1} )\ootimes (\vex_i-\vex_{i+1} ), 
\\
    \pp_q &  \coloneq \ppp_q - \chi_i(1-\chi_i)\Big( |\vex_i - \vex_{i+1}|^2  - \int_{\mathbb T^3} |\vex_i - \vex_{i+1}|^2 \dd x\Big), 
\end{align*}
where we have used the inverse divergence operator $\mathcal R$ from Section \ref{s:inverse-div}. It is easy to check that $\vex_i-\vex_{i+1}$ has mean zero (so that $\RRR_q$ is well-defined), and furthermore, $\pp_q$ has mean zero, while $\RRR_q$ is symmetric and trace-free. A routine computation shows that the functions $(\vv_q, \pp_q, \RRR_q)$ solve:
\begin{equation}
\left\{ \begin{alignedat}{-1}
\del_t \vv_q+\div (\vv_q\otimes \vv_q)  +\nabla \pp_q   &=  \div \RRR_q,
\\
  \nabla \cdot \vv_q &= 0,
  \\
  \vv_q \big|_{t=0}&= \vin * \psi_{\ell_q}.
\end{alignedat}\right.  \label{e:subsol-glued-euler}
\end{equation}
\begin{prop}[Estimates for $\vv_q$]  For  all $N\ge 0$,
\begin{align}\|\vv_{q}-v_{\ell_q}\|_{\alpha} & \lesssim \delta_{q+1}^{1 / 2} \ell_q^\alpha,\label{e:stability-vv_q} \\
\|\vv_{q}-v_{\ell_q}\|_{N+\alpha} & \lesssim \tau_{q} \delta_{q+1} \ell_q^{-1-N+\alpha}, \label{e:stability-vv_q-N} \\
\|\vv_{q}\|_{1+N} & \lesssim \delta_{q}^{1 / 2} \lambda_{q} \ell_q^{-N},\label{e:vv_q-bound}
\end{align}
\end{prop}
\begin{proof}
Since the gluing functions $\chi_i$ do not depend on the space variable,  \eqref{e:stability-v} immediately implies \eqref{e:stability-vv_q-N}. The estimate \eqref{e:stability-vv_q} follows because $\delta_q^{1/2}\tau_q\ell_q^{-1}\le 1$, and \eqref{e:vv_q-bound} follows by the simple bound $\|\vv_q\|_{1+N} \le \|v_{\ell_q}\|_{1+N}+\|v_{\ell_q}-\vv_q\|_{1+N+\alpha}$  since $\ell_q^{-1} \le \lambda_q$.
\end{proof}
\begin{prop}[Estimates for $\RRR_q$]For all $N\ge 0$,\label{p:estimate-RRRq}
    \begin{align}
        \|\RRR_q\|_{N+\alpha}
        &\lesssim \delta_{q+1} \ell_q^{-N+\alpha}, \label{e:RRR_q-N+alpha-bd}
        \\
        \| \matd{\vv_q} \RRR_q\|_{N+\alpha }
        &\lesssim  \delta_{q+1} \delta_{q}^{1/2}\lambda_q \ell_q^{-N-\alpha}.\label{e:matd-RRR_q}
    \end{align}
\end{prop}
\begin{proof}
See \cite[Proposition 4.4]{zbMATH07038033}.
\end{proof}

\begin{prop}[Energy of $\vv_q$] For all $t\in[0,T]$,
\label{p:energy-of-vv_q}
$
   \left| \int_{\mathbb T^3} |\vv_q|^2-|v_{\ell_q}|^2\dd x \right|\lesssim \delta_{q+1}\ell_q^{\alpha},
$
and hence, for $a\gg 1$, and for $t\in[1-\tau_{q-1},T]$,
$
    \frac{\delta_{q+1}}{2\lambda_q^\alpha} \le e(t) - \int_{\mathbb T^3} |\vv_q|^2 \dd x \le 2 \delta_{q+1}.
$
\end{prop}
\begin{proof}
See \cite[Proposition 4.5]{zbMATH07038033}.
\end{proof}

\subsection{Definition of velocity increment $w_{q+1}$}
\label{ss:defn-of-wq+1}
\subsubsection{Space-time cutoffs}
\label{sss:time-cutoffs}
Define the index
\[ i_q \coloneq  \bigg\lfloor \frac1{\tau_q} \bigg\rfloor - 2\in\mathbb Z_{\ge 0}.  \]
 We will define the cutoffs $\eta_i$ by
\begin{align}
    \eta_i (x,t) \coloneq \begin{cases}
 \bar\eta_i(t) & 0\le i< i_q, \\
 \tilde\eta_i(x,t) & \phantom{0\le{}}i\ge i_q,
 \end{cases} \label{e:defn-cutoff-total}
\end{align}
where $\tilde\eta_i$ are `squiggling space-time cutoffs' we will define shortly, and $\bar\eta_i$ are `straight cutoffs' defined as follows: let $\bar\eta_0\in C_c^\infty(J_0\cup I_0 \cup J_1;[0,1])$  satisfy
\[  \supp \bar\eta_0 = I_0 + \Big[  -\frac{\tau_q}6,\ \frac{\tau_q}6\Big] = \Big[  \frac{\tau_q}3-\frac{\tau_q}6,\ \frac{2\tau _q}3+\frac{\tau_q}6\Big]\supset I_0, \]
be identically 1 on $I_0$, and satisfy the derivative estimates for $N\ge 0$:
\[ \|\del_t ^N \bar\eta_0\|_{C^0_t} \lesssim \tau_q^{-N}. \]
Then we set $\bar\eta_i(t)\coloneq\bar\eta_0(t-t_i)$ for $0\le i\le i_q$. Next, we define $\tilde\eta_i$.
\subsubsection{Squiggling space-time cutoffs}
\label{sss:squiggling-cutoffs}
These cutoffs are adapted from \cite[Section 5.2]{zbMATH07038033}. Let $\epsilon \in (0,\frac13)$,  $\epsilon_0\ll 1$ and define for $i_q\le i\le i_{\textup{max}}$ using the mollifiers \eqref{e:defn-mollifier-t}, \eqref{e:defn-mollifier-x},
\begin{align*}
     I_{i}'&\coloneq I_i + \Big[-\frac{(1-\epsilon)\tau_q}3, \frac{(1-\epsilon)\tau_q}3\Big] = \Big[i\tau_q+\frac{\epsilon\tau_q}3, i\tau_q+ \frac{(3-\epsilon)\tau_q}3\Big]  ,\\
         I_i'' &\coloneq \Big\{\Big(x,t+\frac{2\epsilon\tau_q}3\sin (2\pi x_1) \Big) : x\in\mathbb T^3, \ t \in I_i'\Big\}\subset  \mathbb T^3 \times \mathbb R,   \\
         \tilde \eta_i (x,t)&\coloneq \ind_{I_i''}*_x\psi_{\epsilon_0} *_t \phi_{\epsilon_0 \tau_q} = \frac1{\epsilon_0^3}\frac1{\epsilon_0 \tau_q} \iint_{I_i''}\psi\Big(\frac{x-y}{\epsilon_0}\Big)\phi\Big(\frac{t-s}{\epsilon_0\tau_q}\Big)\dd y \dd s.
\end{align*}
\begin{figure}
\begin{center}
    \begin{tikzpicture}

\begin{axis}[clip=false,width=7cm,height=6cm,
    axis line style={draw=none}, 
    disabledatascaling,
    yticklabels={,,},     
    xticklabels={,,},
    tick style={draw=none},
    ]


\addplot [domain=0:1,samples=100, name path global = A](
    {\crviii{x}},
	{x}
	)
	node[%
        pos=0.5%
        ]%
        (curvelabel3)
        { }%
        ;
        
\addplot [domain=0:1,samples=100, name path global = B](
    {\crviv{x}},
	{x}
	)
	node[%
        pos=0.5%
        ]%
        (curvelabel4)
        {}%
        ; 
\addplot [\colorEta] fill between [of=A and B];
\fill[ fill=red!40!white] (\varTAU-\varEPS/3 *\varTAU,0) -- (\varTAU+\varEPS/3 *\varTAU,0.25) -- (\varTAU-\varEPS/3 *\varTAU,0.5) -- cycle; 

\draw (\varTAU-\varEPS*\varTAU,0.75) circle (1.5pt) node [left] {$p$} ; 
\addplot [domain=-0.035:1,samples=2]( 
    {\varTAU},
	{x}
	)  ;
\addplot [domain=-0:1,samples=2, dashed]( 
    {(1-\varEPS/3) * \varTAU},
	{x}
	)  ;
	\addplot [domain=-0:1,samples=2, dashed]( 
    {(1+\varEPS/3) * \varTAU},
	{x}
	)  ;
\addplot [domain=0:((1+\varEPS/3) * \varTAU) ,samples=2, dashed]( 
    {x},
	{0.5}
	)  ;
\addplot [domain=0:1,samples=2, name path global = A](
    {\varTAU/3},
	{x}
	)
        ;
\addplot [domain=0:1,samples=2, name path global = B](
    {2*\varTAU/3},
	{x}
	)
	node[%
        pos=0.5%
        ]%
        (curvelabel)
        { }%
        ;
\addplot [\colorR] fill between [of=A and B];

\draw[->] (-0.5*\varTAU,0) -- (1.4*\varTAU,0) node[right] {$t$};

\node at (0,0) [below]{$t_i$};
\node at (\varTAU,0) [below]{$t_{i+1}$};
\draw [->] (0,-0.035) -- (0,1+\varTAU/2) node [above] {$x_1$} ;
\draw (\varTAU/20,1) -- (-\varTAU/20,1) node [left] {$1$} ;
\draw (\varTAU/20,0.5) -- (-\varTAU/20,0.5) node [left] {$\frac12$} ; 
\draw[|<->|] (\varTAU-\varEPS*\varTAU/3,1+\varTAU/5) -- (\varTAU+\varEPS*\varTAU/3,1+\varTAU/5);
\node at (\varTAU,1+\varTAU/5) [above]{\footnotesize$ 2\epsilon \tau_q/3$};

\draw[<->](-\varTAU/3,-1.5*\varTAU) -- (\varTAU/3,-1.5*\varTAU) node[midway,below] {$J_i$};
\draw[<->](\varTAU-\varTAU/3,-1.5*\varTAU) -- (4*\varTAU/3,-1.5*\varTAU) node[midway,below]  {$J_{i+1}$};

\draw[\colorR, <->](\varTAU/3,-\varTAU/2) -- (\varTAU-\varTAU/3,-\varTAU/2) node[midway,below ]{$I_i$};
\begin{scope}
\clip (\varTAU-\varEPS/3 *\varTAU,0) -- (\varTAU+\varEPS/3 *\varTAU,0.25) -- (\varTAU-\varEPS/3 *\varTAU,0.5) -- cycle; 

\draw[red,very thick] (\varTAU-\varEPS/3 *\varTAU,0) -- (\varTAU+\varEPS/3 *\varTAU,0.25) -- (\varTAU-\varEPS/3 *\varTAU,0.5) -- cycle; 
\end{scope}

\end{axis} 
\end{tikzpicture}
\caption{The set $I''_i$ (light blue $\color{\colorEta}\rule{1ex}{1ex}$). To ensure that $I_i\times\mathbb T^3$ (blue $\color{\colorR}\rule{1ex}{1ex}$) is a subset of $I''_i$, we require the leftmost point of the right boundary, $p=(t_{i+1}-\frac{2\epsilon\tau_q}3,\frac34)$ to be outside $I_i\times\mathbb T^3 $: this is ensured by $\epsilon<1/3$. The triangle (red $\color{red!40!white}\rule{1ex}{1ex}$)
  is a visual proof that any vertical slice of $I''_i $ for $t\in[t_i,t_{i+1}]$ has length at least $1/4$. Point 4 of Lemma \ref{l:squiggle-properties} then follows by taking $\epsilon_0\ll1$.}
\label{f:eta-tilde}
\end{center}
\end{figure}
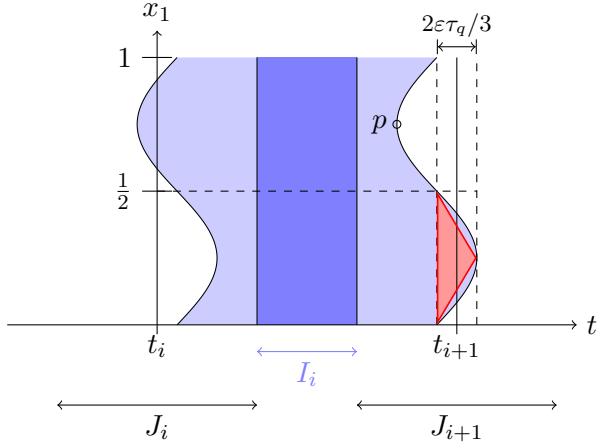

\begin{lem} \label{l:squiggle-properties} The  functions $\eta_i$ satisfy  for all $i=0,1,\dots,i_{\textup{max}}$:
    \begin{enumerate}
    \item $\eta_i \in C^\infty_c(\mathbb T^3 \times (J_i\cup I_i\cup J_{i+1}) ; [0,1])$, with the estimates for $n,m\ge 0$:\begin{align*}
         \|\del^n_t  \eta_i\|_{L^\infty_tC^m_x    }  \lesssim_{n,m} \tau_q^{-n}. 
    \end{align*}
    \item $\eta_i(\cdot,t) \equiv 1$ for $t\in I_i$.
    \item $\supp \eta_i$ are pairwise disjoint.
    \item  For all $t\in[t_{i_q},T]$, $c_\eta \le \sum_{i=i_q}^{i_{\textup{max}}}\int_{\mathbb T^3}  \eta_i^2(x,t) \dd x \le 1$ for a constant $c_\eta$ independent of $q$ (one can take $c_\eta=1/5$).
    \item For all $0\le     i<i_q$, $\eta_i $ does not depend on $x$, and  is supported in time on the $\tfrac{\tau_q}6$-neighbourhood of $I_i$.
\end{enumerate}%
\end{lem}%
\subsubsection{Energy gap decomposition}
\label{ss:energy-gap}

Let $\zeta$ be a smooth cutoff function such that
\begin{align*}
    t\le t_{i_q}&\implies \zeta\equiv 1,  
    \\
    t\ge t_{i_q}+\tau_q/3 &\implies \zeta\equiv 0,\\
    t\in (t_{i_q},t_{i_q}+\tau_q/3) &\implies |\del_t^N\zeta| \lesssim_N \tau_q^{-N}.
\end{align*}
Observe that $t_{i_q}=\tau_q \lfloor \tau_q^{-1}\rfloor-2\tau_q\in (1-3\tau_q,1-2\tau_q)$.
This function is needed for because we do not have the estimate \eqref{e:energy-q-estimate} for small times. For $a\gg1$, $1-3\tau_q > 1-\tau_{q-1}$, so we define the energy gap $\rho_q$ as
\[  \rho_q(t) \coloneq  \frac{\delta_{q+1}}2\zeta(t) + \frac{1-\zeta(t)}3\Big(e(t) - \frac{\delta_{q+2}}2 - \int_{\mathbb T^3} |\vv_q|^2 \dd x \Big). \]Note that
\begin{align}
  t\in[t_{i_q+1},T]\implies   \rho_q(t) &=  \frac13\Big(e(t) - \frac{\delta_{q+2}}2 - \int_{\mathbb T^3} |\vv_q(x,t)|^2 \dd x \Big). \label{e:energy-gap}
\end{align}
We next use the cutoff functions $\zeta$ and $\eta_i$ to split $\rho_q$ into functions $\rho_{q,i}$ supported on disjoint regions of space-time.
Note $t_{i_q+1}=\tau_q (\lfloor \tau_q^{-1}\rfloor -1) \le 1-\tau_q$, so that $t\ge 1-\tau_q$ implies $\zeta(t)=0$.
We decompose $\rho_q$ by setting ($i=0,1,\dots,i_{\textup{max}}$)
\begin{align}
 \rho_{q,i} (x,t) &\coloneq \frac{\eta^2_i(x,t)}{ \zeta(t) + \sum_{i=i_q}^{i_{\textup{max}}} \int_{\mathbb T^3}\eta_i^2 (\tilde x,t)\dd \tilde x }\rho_q(t).  \label{e:energy-q-i-gap}
\end{align}
Note that for $t\in I_i$, $i\ge 0$,
\[ \zeta(t) + \sum_{i=i_q}^{i_{\textup{max}}} \int_{\mathbb T^3} \eta_i^2 (\tilde x,t)\dd \tilde x \equiv 1, \]
and for $t\in J_i$ there is  the strictly positive lower bound
\[ \zeta(t) + \sum_{i=i_q}^{i_{\textup{max}}} \int_{\mathbb T^3} \eta_i^2 (\tilde x,t)\dd \tilde x   \ge c_\eta . \] In particular, $\rho_{q,i}$ is well-defined even for times $t$ where $ \sum_{i= 0}^{i_{\textup{max}}} \int_{\mathbb T^3}\eta_i^2 (\tilde x,t)\dd \tilde x=0$. By construction, $\rho_q = \int_{\mathbb T^3}\sum_{i= 0}^{i_{\textup{max}}} \rho_{q,i} \dd x $ for times $t\ge1-\tau_q$.
\begin{figure}
    \begin{center}
        \begin{tikzpicture}

\begin{axis}[clip=false,width=13cm,height=7.5cm,
    axis line style={draw=none}, 
    disabledatascaling,
    yticklabels={,,},     
    xticklabels={,,},
    tick style={draw=none},
    ]
 \addplot [line width = 0.4cm,domain=0:2*pi,samples=100](
	{\crvi{x}},
	{\crvii{x}}
	); 
 \addplot [white, line width = 0.35cm,domain=0:2*pi,samples=100](
	{\crvi{x}},
	{\crvii{x}}
	); 
 
\foreach \k [count=\K] in {0,...,\varN}
{
\ifnum \k<\varNi 
\addplot [domain=0:1,samples=100, name path global =A](%
    {\k*\varTAU+\crvi{x}},%
	{x}%
	)%
	node[%
        pos=0.5%
        ]%
        { }%
        ;
\addplot [domain=0:1,samples=100, name path global =B](
    {\k*\varTAU+\crvii{x}},
	{x}
	)
	node[%
        pos=0.5%
        ]%
        {}%
        ; 
\addplot [\colorEta] fill between [of=A and B];
\else 
\addplot [domain=0:1,samples=100, name path global = A](
    {\k*\varTAU+\crviii{x}},
	{x}
	)
	node[%
        pos=0.5%
        ]%
        { }%
        ;
        
\addplot [domain=0:1,samples=100, name path global = B](
    {\k*\varTAU+\crviv{x}},
	{x}
	)
	node[%
        pos=0.5%
        ]%
        {}%
        ; 
\addplot [\colorEta] fill between [of=A and B];
\fi 
\addplot [domain=-0.035:1,samples=2, dashed]( 
    {\k*\varTAU},
	{x}
	)   ;
\addplot [domain=0:1,samples=2, name path global = A](
    {\k*\varTAU+\varTAU/3},
	{x}
	)
        ;
\addplot [domain=0:1,samples=2, name path global = B](
    {\k*\varTAU+2*\varTAU/3},
	{x}
	)
	node[%
        pos=0.5%
        ]%
        { }%
        ;
\addplot [\colorR] fill between [of=A and B];
}
\addplot [domain=-0.035:1,samples=2, dashed]( 
    {\varN*\varTAU+\varTAU},
	{x}
	)
	node[%
        pos=0.5%
        ]%
        { }%
        ;

\draw[->] (-0.15*\varTAU,0) -- (\varN*\varTAU+1.4*\varTAU,0) node[right] {$t$};
\node at (\varNi*\varTAU,0) [below]{$t_{i_q}$};
\node at (\varNi*\varTAU+ \varTAU,0) [below]{$t_{i_q+1}$};
\node at (\varNi*\varTAU+2*\varTAU,0) [below]{$\tau_q\lfloor\tau_q^{-1}\rfloor$};

\draw [->] (0,0) -- (0,1+\varTAU/2) node [above] {$x_1$} ;
\draw (\varTAU/20,1) -- (-\varTAU/20,1) node [left] {$1$} ;

\node at (\varNi*\varTAU+\varTAU/3+\varTAU/6-\varTAU,1+\varTAU) (a) [above, \colorR]{$I_i$};
\draw [\colorR,->] (a) -- (\varTAU*2.5-\varTAU,1+\varTAU/6);
\draw [\colorR,->] (a) -- (\varTAU*3.5-\varTAU,1+\varTAU/10);
\draw [\colorR,->] (a) -- (\varTAU*4.5-\varTAU,1+\varTAU/6); 

\node at (\varNi*\varTAU+\varTAU*2-\varTAU,1+\varTAU ) (b) [above]{$J_i$};
\draw [thick, white] (b) -- (\varNi*\varTAU,1+\varTAU/10) node (c) {} ;
\draw [] (b) -- (\varNi*\varTAU,1+\varTAU/10) node (c) {} ;
\draw [<->] ($(c)-(\varTAU/3,0)$) -- ($(c)+(\varTAU/3,0)$);
\draw [thick, white] (b) -- (\varNi*\varTAU+2*\varTAU-\varTAU,1+\varTAU/10) node (c) {} ;
\draw [] (b) -- (\varNi*\varTAU+2*\varTAU-\varTAU,1+\varTAU/10) node (c) {} ;
\draw [<->] ($(c)-(\varTAU/3,0)$) -- ($(c)+(\varTAU/3,0)$);
\draw [thick, white] (b) -- (\varNi*\varTAU+3*\varTAU-\varTAU,1+\varTAU/10) node (c) {} ;
\draw [] (b) -- (\varNi*\varTAU+3*\varTAU-\varTAU,1+\varTAU/10) node (c) {} ;
\draw [<->] ($(c)-(\varTAU/3,0)$) -- ($(c)+(\varTAU/3,0)$);


\draw [\colorEtaMinusOne, very thick] (0,\varYOFF+1) -- (\varETASTART +0.0005,\varYOFF+1);
\addplot [\colorEtaMinusOne, very thick, domain=\varETASTART :(\varETASTART +\varTAU/3),samples=100] ( 
    {x},
	{\etaminusone{(x-(\varETASTART) )/(\varTAU/3) } +\varYOFF}
	)
	node[%
        pos=0.5%
        ]%
        { }%
        ;
\draw [\colorEtaMinusOne, very thick] (\varETASTART +\varTAU/3-0.0005,\varYOFF) -- (\varN*\varTAU+\varTAU,\varYOFF)  ;
\draw[->] (-0.15*\varTAU,\varYOFF ) -- (\varN*\varTAU+1.4*\varTAU,\varYOFF ) node[right] {$t$};
\node at (\varNi*\varTAU, \varYOFF ) [below]{$t_{i_q}$};
\node at (\varNi*\varTAU+\varTAU,\varYOFF ) [below]{$t_{i_q+1}$};
\node at (\varNi*\varTAU+2*\varTAU,\varYOFF ) [below]{$\tau_q\lfloor\tau_q^{-1}\rfloor$};

\draw [->] (0,\varYOFF -0.05) -- (0,\varYOFF + 1+\varTAU/2 ) node [above] {$\zeta(t)$} ;

\addplot [domain=-0.035:1-0.2*\varYOFF,samples=2, dashed]( 
    {\varNi*\varTAU},
	{x+\varYOFF}
	)
	node[%
        pos=0.5%
        ]%
        { }%
        ;
\addplot [domain=-0.035:1-1*\varYOFF,samples=2, dashed]( 
    {\varNi*\varTAU+\varTAU/3},
	{x+\varYOFF}
	)
	node[%
        pos=0.5%
        ]%
        { }%
        ;

\addplot [domain=-0.035:1-0.2*\varYOFF ,samples=2, dashed]( 
    {\varNi*\varTAU+\varTAU},
	{x+\varYOFF}
	)
	node[%
        pos=0.5%
        ]%
        { }%
        ; 
\addplot [domain=-0.035:1-0.2*\varYOFF ,samples=2, dashed]( 
    {\varNi*\varTAU+2*\varTAU},
	{x+\varYOFF}
	)
	node[%
        pos=0.5%
        ]%
        { }%
        ;
\draw (\varTAU/20,1+\varYOFF) -- (-\varTAU/20,1+\varYOFF) node [left] {$1$} ;
\end{axis} 
\end{tikzpicture}
    \caption{At time $t_{i_q}$, we switch from using straight cutoffs $\overline \eta_i$ to the squiggling cutoffs $\tilde\eta_i$. We give ourselves a little time to `turn off' $\zeta$ since we no longer need it to ensure $\rho_{i,q}$ is well-defined. At time $t_{i_q+1}=\tau_q(\lfloor \tau_q^{-1}\rfloor-1)< 1$, we begin controlling the energy. }
    \label{f:all-eta}
    \end{center}
\end{figure}
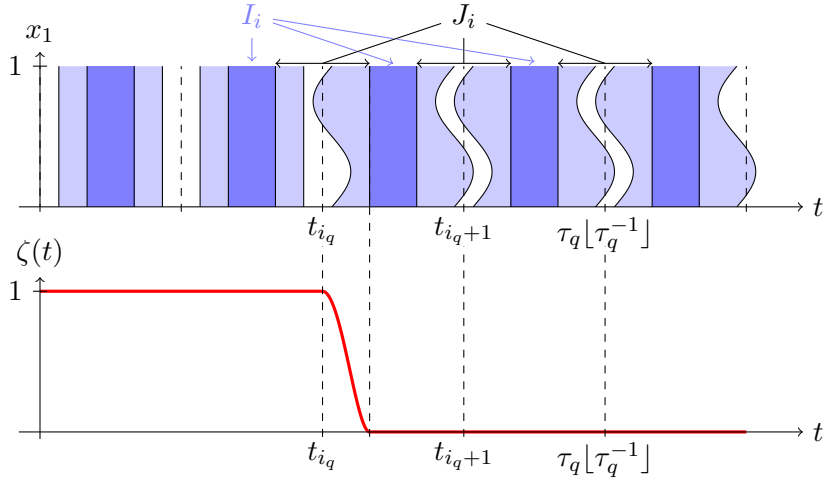

\begin{prop}[Estimates for $\rho_q$ and $\rho_{q,i}$]\label{p:estimates-for-energy-gap} For all  $t \in[0,T]$,
    \begin{align} \frac{\delta_{q+1}}{8 \lambda_{q}^{\alpha}} \leq \rho_{q}(t) & \leq \delta_{q+1}  \label{e:rho_q-upper-lower-bd}, \\
    \|\rho_{q, i}\|_{0} & \leq {\delta_{q+1}}/{c_{\eta}} ,\label{e:rho_qi-C0-bd} \\
    \|\rho_{q, i}\|_{N} & \lesssim \delta_{q+1} \label{e:rho_qi-CN-bd}, \\
    \|\rho_{q, i}^{1/2}\|_{N} & \lesssim \delta_{q+1}^{1/2} \label{e:rho_qi-sqrt-CN-bd}, \\
    \sup_{t\in[0,T]}|\partial_{t} \rho_{q}(t)| & \lesssim \delta_{q+1} \tau_q^{-1} \label{e:rho_q-ddt-C0},\\
    \sup_{t\in \cup_{i\ge0} I_i} |\partial_{t} \rho_{q}(t)| & \lesssim \delta_{q+1} \delta_q^{1/2}\lambda_q \label{e:rho_q-ddt-C0-better-on-Ii},\\
    \|\partial_{t} \rho_{q, i}\|_{N} & \lesssim \delta_{q+1} \tau_{q}^{-1} \label{e:rho_qi-ddt-CN}. \end{align}
\end{prop}
\begin{proof} For $t\ge1-\tau_q$,
\eqref{e:rho_q-upper-lower-bd} follows from because of $4\delta_{q+2}\le \delta_{q+1}$, which holds for $a\gg1$. Then, \eqref{e:rho_qi-C0-bd} follows from point 4 of Lemma \ref{l:squiggle-properties} and the definition of $\zeta$. The implicit constant in \eqref{e:rho_qi-CN-bd} depends only on $\epsilon_0$ which is a fixed universal constant, and not on $\tau_q$ due to the definition of $\tilde\eta_i$. For \eqref{e:rho_qi-sqrt-CN-bd} we emphasise that this norm is $L^\infty$ in time and $C^N$ in space. The remaining estimates follow from the definition \eqref{e:tau_q-and-t_i} of $\tau_q$.

For $t<1-\tau_q$, all of the estimates are even more straightforward. We stress that  \eqref{e:rho_qi-C0-bd} still holds, precisely because we added the term $\zeta$ in the denominator.
\end{proof}
We remark that the analogous estimate of \eqref{e:rho_q-ddt-C0} is worse than \eqref{e:rho_q-ddt-C0-better-on-Ii} (which is \cite[(5.14)]{zbMATH07038033}) by a small $\alpha$ power. This estimate is only used to show \eqref{e:tilde-Rqi-N-matd}. We kept \eqref{e:tilde-Rqi-N-matd} the same as the corresponding estimate  \cite[(5.38)]{zbMATH07038033} for the convenience of the reader.

\subsubsection{Inverse flow map}
The inverse flow (also called `back-to-labels') maps $\Phi_i$ are defined as the solutions to the vector transport equation \label{ss:back-to-labels}
    \begin{align}
         (\partial_t+\vv_q \cdot \nabla) \Phi_i = 0, \quad \Phi_i\big|_{t=t_i} = x. \label{e:phi_i-defn}
    \end{align}
For each fixed $t$, $\Phi_i$ is the inverse mapping for the flow   $\Xi=\Xi(x,t)$, i.e. the solution to $\del_t \Xi (x,t) = \vv_q (\Xi(x,t),t) $ with initial data $\Xi\big|_{t=t_i} = x$.

In preparation for the use of the geometric property \eqref{e:mikado-geometric-property} in Lemma \ref{l:mikado}, we define
\[
    R_{q,i} \coloneq I_{3\times 3} - \frac{\eta_{i}^2 \RRR_q}{\rho_{q,i}}
\]
In the above, we multiply by $\eta_{i}^2$ to select a component of $\supp\RRR_q$; note that $\eta_i^2\RRR_q$ is compactly supported in $I_i$, where it is identically equal to $\RRR_q$.

In order to use  $C_k\cdot k=0$ from  \eqref{e:mikado-C_k},  we have to conjugate with $\nabla \Phi$:
\begin{align} \widetilde R_{q,i}\coloneq \nabla\Phi R_{q,i}\nabla \Phi^\TT  =\nabla\Phi \nabla \Phi^\TT  - \frac1{\rho_{q,i}}\eta_i^2 \nabla\Phi \RRR_{q}\nabla \Phi^\TT . \label{e:defn-Rqi} \end{align}
See the algebraic computation in \eqref{e:cancellation}.

\begin{prop}[Estimates for $\Phi_i$] For $a\gg 1$,
\label{p:estimates-for-inverse-flow-map}\begin{align}
 t\in J_i=\supp \eta_i   \implies     \|\nabla\Phi_i-I_{3\times3}\|_0 &\le\frac1{10}.\label{e:nabla-phi-i-minus-I3x3}
 \\
 \| (\nabla \Phi_i)^{-1}\|_N+ \|  \nabla \Phi_i\|_N &\le \ell_q^{-N}, \label{e:nabla-phi-i-CN}
\\
\|\matd{\vv_q} \nabla \Phi_i\|_N &\lesssim \delta_q^{1/2} \lambda_q \ell_q^{-N}. \label{e:nabla-phi-i-matd}
\end{align}
\end{prop}
\begin{proof}
The first two follow because there exists $C=C_N$ such that
\begin{align*}
    \|\nabla \Phi_i - I_{3\times 3}\|_0
    &\leq \ee^{(t-t_i)\|\vv_q\|_{1}}-1 \le \tau_q \|\vv_q\|_{1} \lesssim \ell_q^{2\alpha} \ll 1,
    \\
[ \Phi_i - I_{3\times 3}]_N &\leq C (t-t_i)[\vv_q]_N\ee ^{C(t-t_i)[\vv_q]_1},
\end{align*}
 as in \cite[App. D]{zbMATH06456007}. $a\gg1$ is used to control constants.

 For the material derivative, we compute \begin{align*}
    (\matd{\vv_q} \nabla \Phi_i)_{ab}
    &= \matd{\vv_q} \del_{b} (\Phi_i)_a
    =\del_t \del_{b} (\Phi_i)_a  +(\vv_q)_c\del_{c} \del_{b} (\Phi_i)_a\\
    &=\del_{b}(\del_t  (\Phi_i)_a  +(\vv_q)_c\del_{c}  (\Phi_i)_a) - \del_{b}(\vv_q)_c \del_{c}  (\Phi_i)_a \\
    &=  \del_{b} ( \matd{\vv_q} (\Phi_i)_a) - (\nabla \vv_q^\TT \nabla \Phi_i)_{ab},
\end{align*}
and hence by \eqref{e:phi_i-defn}, $\matd{\vv_q} \nabla \Phi_i= -\nabla \vv_q^\TT \nabla \Phi_i$.
So \eqref{e:nabla-phi-i-matd} follows from \eqref{e:nabla-phi-i-CN} and \eqref{e:vv_q-bound}.
\end{proof}

\begin{prop}[Estimates for $\widetilde R_{q,i}$] For $a\gg1$,
\begin{align}
    \|\widetilde R_{q,i}-I_{3\times3}\|_0&\le \frac12,\label{e:tilde-Rqi-minus-identity}
    \\
    \|\widetilde R_{q,i}\|_N&\lesssim \ell_q^{-N}, \label{e:tilde-Rqi-N}
    \\
    \|\matd{\vv_q} \widetilde R_{q,i}\|_N &\lesssim \tau_q^{-1} \ell_q^{-N}. \label{e:tilde-Rqi-N-matd}
\end{align}
\end{prop}
\begin{proof}
For the first two estimates, we write
\begin{align*}
    &\widetilde R_{q,i}-I_{3\times 3}
   =\nabla \Phi_i (R_{q,i} - I_{3\times 3} ) \nabla\Phi^\TT + \nabla\Phi \nabla\Phi^\TT - I_{3\times 3}
    \\
    &= \nabla \Phi \frac{\RRR_q}{\rho_{q,i}}\nabla\Phi^\TT \! + (\nabla\Phi-I_{3\times 3})(\nabla\Phi^\TT\! +I_{3\times 3}) + (\nabla\Phi-I_{3\times 3}) -(\nabla\Phi^\TT\!-I_{3\times 3}),
\end{align*}
so that when $a\gg 1$ so that $ \lambda_q^\alpha \ell_q^\alpha < \frac{9}{100}$, by \eqref{e:RRR_q-N+alpha-bd} and the fact that either $\RRR_q\equiv0$ or $\rho_{q,i}=\rho_q>\lambda_q^{-\alpha}\delta_{q+1}$ by \eqref{e:rho_q-upper-lower-bd}, \begin{align*}
    \|\widetilde R_{q,i}-I_{3\times 3}\|_{C^0} \le \lambda_q^{\alpha} \ell_q^{\alpha} + \frac1{10}\Big(2+\frac1{10}\Big) + \frac1{10}+\frac1{10} < \frac12.
\end{align*}
The second estimate \eqref{e:tilde-Rqi-minus-identity} follows by routine direct computation. For the estimate on the material derivative, as $\eta^2_i\equiv 1$ and $\rho_{q,i}=\rho_q$ on a neighbourhood of the support of $R_{q,i}$, we have
\begin{align*}
     \matd{\vv_q} R_{q,i}  &= \matd{\vv_q} (\rho_{q}^{-1}\RRR_q)
    = -\frac{\del_t \rho_{q}}{\rho_{q}^2} \RRR_q + \rho_{q} \matd{\vv_q} \RRR_q,
\end{align*}
Hence, from the estimates \eqref{e:matd-RRR_q}, \eqref{e:rho_q-upper-lower-bd}, and \eqref{e:rho_q-ddt-C0-better-on-Ii} (since $\RR_q$ is supported on $\mathbb T^3\times\bigcup_{i\ge0}I_i$), we obtain
\begin{align*}
    \|\matd{\vv_q} R_{q,i}\|_N \le \tau_q^{-1}\ell_q^{-N}.
\end{align*}
As $\|\nabla \Phi_i\|_0\lesssim 1$ by \eqref{e:nabla-phi-i-minus-I3x3}, the estimate \eqref{e:nabla-phi-i-matd} follows:
\begin{align*}
   & \| \matd{\vv_q} \widetilde R_{q,i}\|_N \\
   &\lesssim \| \matd{\vv_q}\!\nabla \Phi^\TT  R_{q,i} \nabla\Phi\|_N +\| \nabla \Phi  \matd{\vv_q}R_{q,i}  \nabla\Phi^\TT \|_N + \| \nabla \Phi  R_{q,i} \matd{\vv_q} \!\nabla\Phi^\TT \|_N
    \\
    &\lesssim \| \matd{\vv_q} \nabla \Phi \|_N \| R_{q,i}\|_0 + \| \matd{\vv_q} \nabla \Phi \|_0 \| R_{q,i}\|_N 
    \\
    &\quad + \| \matd{\vv_q} \nabla \Phi \|_0 \| R_{q,i}\|_0 \|\nabla\Phi\|_N + \|\nabla\Phi\|_N \| \matd{\vv_q} R_{q,i}\|_0 +  \| \matd{\vv_q} R_{q,i}\|_N
    \\
    &\lesssim \tau_q^{-1} \ell_q^{-N}. \qedhere
\end{align*}
\end{proof}

\subsection{Mikado flows}
\label{ss:mikado}

In this subsection, we summarize the Mikado flows as presented in \cite{zbMATH07038033} as a lemma. Mikado flows were first introduced in \cite{zbMATH06710292}. Write $\mathcal S^{3\times 3}_+$ for the space of symmetric positive definite $3\times 3$ matrices, and $\overline{\mathbb B}_{1/2} (I_{3\times 3}) \subset \mathcal S^{3\times 3}_+$ for the closed ball of radius $1/2$ around the identity matrix.
\begin{lem}
\label{l:mikado}
There exists a smooth map
\[
W: \overline{\mathbb B}_{1/2} (I_{3\times 3}) \times \mathbb{T}^{3} \rightarrow \mathbb{R}^{3}
\]
such that, for every fixed $R \in \overline{\mathbb B}_{1/2} (I_{3\times 3})$, $W(R,\cdot)$ is a presureless solution to the stationary incompressible Euler equations,
\begin{align*}
\operatorname{div}_{\xi}(W(R, \xi) \otimes W(R, \xi))&=0,  
\\
\operatorname{div}_{\xi} W(R, \xi)&=0,  
\end{align*}
and the low frequencies of $W$ are such that
\begin{align}
    \dashint_{\mathbb{T}^{3}} W(R, \xi) \dd \xi &=0, 
     \\
\dashint_{\mathbb{T}^{3}} W(R, \xi) \otimes W(R, \xi) \dd \xi &=R .  \label{e:mikado-geometric-property}
\end{align}
It follows from the smoothness of $W$ that there exist  $a_k\in C^\infty(\overline{\mathbb B}_{1/2} (I_{3\times 3});\mathbb R^3)$ such that
\begin{align}
    W(R, \xi)=\sum_{k \in \mathbb{Z}^{3} \backslash\{0\}} a_{k}(R) \ee^{\ii k \cdot \xi}, \quad a_k(R)\cdot k=0 ,  \label{e:mikado-a_k}
\end{align}
and there exist  $C_k\in C^\infty(\overline{\mathbb B}_{1/2} (I_{3\times 3});\mathbb R^{3\times 3})$ such that \begin{align}
    W(R, \xi)\otimes W(R, \xi)=R +\sum_{k \in \mathbb{Z}^{3} \backslash\{0\}} C_{k}(R) \ee^{\ii k \cdot \xi},\quad  C_k(R) k=0. \label{e:mikado-C_k}
\end{align}
For all $N,m\ge 0$ and $k\in\mathbb Z^3\setminus \{0\}$, there exist constants $C_{N,m}$ such that:
\begin{align}
    \sup _{R \in \overline{\mathbb B}_{1/2} (I_{3\times 3})}\left|D_{R}^{N} a_{k}(R)\right| + \sup _{R \in \overline{\mathbb B}_{1/2} (I_{3\times 3})}\left|D_{R}^{N} C_{k}(R)\right| \le \frac{C_{N,m}}{|k|^{m}},
\label{e:mikado-a_k and c_k estimates}
\end{align}
Here, $D_R$ denotes the derivative with respect to $R$.
\end{lem}
\subsubsection{Principal part and incompressibility corrector}
Define the \emph{principal part} of the perturbation
\label{ss:principal-part-and-incomp-corrector}
\begin{align*}
    \wpq \coloneq \sum_{i= 0}^{i_{\textup{max}}} \wpq[q+1,i]  \coloneq \sum_{i= 0}^{i_{\textup{max}}} \rho_{q,i}^{1/2} \nabla \Phi_i^{-1} W(\widetilde R_{q,i} , \ \lambda_{q+1} \Phi_i).
\end{align*}
We need to define a small corrector $\wcq$ so that $w_{q+1}=\wpq+\wcq$ is divergence-free, i.e. that $w_{q+1}$ is a curl.
 We start by using  the  Fourier decomposition \eqref{e:mikado-a_k} to write
\[ \wpq = \sum_{k,i} \rho_{q,i}^{1/2} \nabla \Phi_i^{-1} a_k(\widetilde R_{q,i} ) \ee^{\ii \lambda_{q+1} k \cdot \Phi_i},\]
where the sum is over $k\in\mathbb Z^3\setminus\{0\}$ and $i\in\{0,1,\dots,i_{\textup{max}}\}$. From $a_k\cdot k=0$ in \eqref{e:mikado-a_k}, it follows that $a_{k}
= -  \frac{ k \times a_{k}}{|k|^{2}} \times  k  $.  Also recall the standard cross product identity\footnote{This easily follows from  $\det M = \epsilon_{ijk} M_{1 i} M_{2j}M_{3k}$, $(\operatorname{cof}M)_{ij}=\frac12 \epsilon_{ipq}\epsilon_{jkl} M_{pk} M_{ql}$ and $A^{-1}=\frac1{\det A} \operatorname{cof} A^\TT $.  }
$M^{-1}(a\times b)=  (\det M) (M^\TT a)\times (M^\TT  b).$ As $\det \nabla \Phi_i=1$, it follows that
\begin{align*}
     \wpq
     &= \frac{-\ii}{\lambda_{q+1}}\sum_{k,i} \rho_{q,i}^{1/2} \frac{\nabla \Phi_i^\TT (k\times a_k(\widetilde R_{q,i} )) }{|k|^2} \times (\nabla  \ee^{\ii \lambda_{q+1} k \cdot \Phi_i}).
\end{align*}

The above considerations show that we can define a divergence-free $v_{q+1}$ with zero mean by setting $ w_{q+1} \coloneq \wpq +\wcq $ and $ v_{q+1} \coloneq \vv_q + w_{q+1}$, where
\begin{align*}
    \wcq &\coloneq  \frac{-\ii}{\lambda_{q+1}}\sum_{k,i}\curl \left( \rho_{q,i}^{1/2} \frac{\nabla \Phi_i^\TT (k\times a_k(\widetilde R_{q,i} )) }{|k|^2} \right)  \ee^{\ii \lambda_{q+1} k \cdot \Phi_i}, 
\end{align*}
Note that $w_{q+1}$ is supported away from zero. Hence,
\[ v_{q+1}\big|_{t=0} = \vv_q\big|_{t=0}=\vin * \psi_{\ell_q}.\]

We define $b_{i,k}=b_{i,k}(x,t)$ and $c_{i,k}=c_{i,k}(x,t)$ by
\begin{align}
    b_{i,k}(x,t)&\coloneq    \rho_{q,i}^{1/2} a_k(\widetilde R_{q,i} ),\label{e:bik-defn} \\
    c_{i,k}(x,t)&\coloneq \frac{-\ii }{\lambda_{q+1}}\curl \left( \rho_{q,i}^{1/2} \frac{\nabla \Phi_i^\TT (k\times a_k(\widetilde R_{q,i} )) }{|k|^2} \right)  \label{e:cik-defn},
\end{align}
so that
\begin{align}
     \wpq&=\sum_{i= 0}^{i_{\textup{max}}} \wpq[q+1,i], & \wpq[q+1,i]&=\sum_{k\neq 0} (\nabla \Phi_i)^{-1}b_{i,k} \ee ^{\ii \lambda_{q+1} k\cdot\Phi_i}, \label{e:wpqi-bik} \\
    \wcq&=\sum_{i= 0}^{i_{\textup{max}}} \wcq[q+1,i], & \wcq[q+1,i]&=\sum_{k\neq 0} c_{i,k} \ee ^{\ii \lambda_{q+1} k\cdot\Phi_i}. \label{e:wcqi-cik}
 \end{align}
\begin{prop}[$C^0$ Estimate for $a_k$ and $b_{i,k}$] There is a universal constant $\overline M$ such that\label{p:b-i-k-C0estimate} 
\begin{align}
    \|b_{i,k}\|_0 \le \frac{\overline M}{|k|^5} \delta_{q+1}^{1/2}. \label{e:b-i-k-C0estimate}
\end{align}
\end{prop}
\begin{proof}
Putting $N=0$ and $m=5$ in \eqref{e:mikado-a_k and c_k estimates} gives  $  \|a_k(\widetilde R_{q,i})\|_0 \le \frac{C_{0,5}}{|k|^5}. $
From the definition of $b_{i,k}$ and the upper bound \eqref{e:rho_qi-C0-bd} $\rho_{q,i} \le  \delta_{q+1}/c_\eta$, we can set $\overline M \coloneq C_{0,5}c_\eta^{-1/2}$.
\end{proof}
We set the constant $M$ in the estimate \eqref{e:vq-C1} as
\begin{align}
M\coloneq 64 \overline M \sum_{k\in\mathbb Z^3\setminus \{0\}} \frac1{|k|^4}<\infty .     \label{d:M}
\end{align}

\begin{prop}[Estimates for $b_{i,k}$ and $c_{i,k}$]
\label{p:estimates-for-b-c}
For $t\in J_i\cup I_i \cup J_{i+1}$,\begin{align}
  \|b_{i,k}\|_N &\lesssim \delta_{q+1}^{1/2} |k|^{-6} \ell_q^{-N} \label{e:bik-CN},
\\
 \|c_{i,k}\|_N &\lesssim \delta_{q+1}^{1/2} |k|^{-6} \ell_q^{-N-1}\lambda_{q+1}^{-1} ,\label{e:cik-CN}
\\
\|\matd{\vv_q} c_{i,k}\|_N &\lesssim \delta_{q+1}^{1/2}|k|^{-6}\ell_q^{-N-1}\lambda_{q+1}^{-1}\tau_q^{-1}.    \label{e:cik-matd-CN}
\end{align}
\end{prop}
\begin{proof}
Similarly to Proposition \ref{p:b-i-k-C0estimate}, the first estimate \eqref{e:bik-CN} follows from the definition \eqref{e:bik-defn} and the estimates \eqref{e:mikado-a_k and c_k estimates}, \eqref{e:rho_qi-sqrt-CN-bd}, \eqref{e:nabla-phi-i-CN}, and  \eqref{e:tilde-Rqi-N}:
\begin{align*}
    \|b_{i,k} \|_N &\lesssim \|\rho_{q,i}^{1/2}\|_N\|a_k(\widetilde R_{q,i})\|_0 +  \|\rho_{q,i}^{1/2}\|_0\|a_k(\widetilde R_{q,i})\|_N
    \\
    &\lesssim \delta_{q+1}^{1/2}\|a_k(\widetilde R_{q,i})\|_0 +  \delta_{q+1}^{1/2}\|a_k(\widetilde R_{q,i})\|_N \lesssim \delta_{q+1}^{1/2}\frac{C_{N,6}}{|k|^6} \ell_q^{-N}.
\end{align*}
A similar calculation gives \eqref{e:cik-CN} from the definition \eqref{e:cik-defn},
\begin{align*}
    \|c_{i,k}\|_N \lesssim \delta_{q+1}^{1/2}\lambda_{q+1}^{-1} \Big\|\eta_i \nabla \Phi_i^T \frac{k\times a_k(\widetilde R_{q,i})}{|k|^2} \Big\|_{N+1}
    \lesssim \lambda_{q+1}^{-1}\delta_{q+1} \ell_q^{-N-1} \frac{C_{N+1,5}}{|k|^6} .
\end{align*}

For \eqref{e:cik-matd-CN}, writing $c_{i,k} = \nabla\times \tilde c_{i,k}$, we commute the curl and compute using \eqref{e:rho_q-ddt-C0}, \eqref{e:nabla-phi-i-matd}, \eqref{e:tilde-Rqi-N-matd} and \eqref{e:mikado-a_k and c_k estimates}:
\begin{align*}
&    \|\matd{\vv_q} c_{i,k} \|_N
\lesssim \|\matd{\vv_q} \tilde c_{i,k}\|_{N+1} + \| \vv_q\|_1 \|\nabla \tilde c_{i,k}\|_{N} + \| \vv_q\|_{N+1} \|\nabla \tilde c_{i,k}\|_{0}
\\
& \lesssim \lambda_{q+1}^{-1} \Bigg \| \frac{\del_t \rho_{q,i}}{\rho_{q,i}^{1/2}}  \nabla \Phi_i^\TT  \frac{  a_k(\widetilde R_{q,i} ) }{|k|} \Bigg \|_N
+ \lambda_{q+1}^{-1}\left\| \rho_{q,i}^{1/2} \matd{\vv_q} \nabla \Phi_i^\TT  \frac{  a_k(\widetilde R_{q,i} ) }{|k|} \right\|_N
\\
&\quad  +\lambda_{q+1}^{-1} \left\| \rho_{q,i}^{1/2} \nabla \Phi_i^\TT  \frac{ D_R a_k(\widetilde R_{q,i} ) }{|k|} \matd{\vv_q}\widetilde R_{q,i} \right\|_N
+ \delta_{q+1}\lambda_q  |k|^{-6} \ell_q^{-N-1}\lambda_{q+1}^{-1}
\\&\lesssim \delta_{q+1}^{1/2} \tau_q^{-1} \lambda_q^{-\alpha}\ell_q^{-N} \lambda_{q+1}^{-1} |k|^{-6}  + \delta_{q+1}\lambda_q\ell_q^{-N}\lambda_{q+1}^{-1} |k|^{-6}
\\&\quad +  \delta_{q+1}^{1/2}\tau_q^{-1}\ell_q^{-N}\lambda_{q+1}^{-1} |k|^{-6}  +  \delta_{q+1}^{1/2} \tau_q^{-1} |k|^{-6} \ell_q^{-N-1}\lambda_{q+1}^{-1}
\\&\lesssim \delta_{q+1}^{1/2}\tau_q^{-1}\ell_q^{-N-1}\lambda_{q+1}^{-1}|k|^{-6}. \qedhere
\end{align*}
\end{proof}
\begin{cor}[Estimates for $\wpq$ and $\wcq$]
\label{c:estimates-for-wpq-wcq}
\begin{align}
    \|\wpq\|_0+\frac1{\lambda_{q+1}} \|\wpq\|_1 &\le \frac M4 \delta_{q+1}^{1/2},\label{e:wpq-decay}
    \\
        \|\wcq\|_0+\frac1{\lambda_{q+1}} \|\wcq\|_1 &\le \ell_q^{-1}\lambda_{q+1}^{-1} \delta_{q+1}^{1/2},\label{e:wcq-decay}
        \\
            \|w_{q+1}\|_0+\frac1{\lambda_{q+1}} \|w_{q+1}\|_1 &\le \frac M2 \delta_{q+1}^{1/2}.\label{e:wq+1-decay}
\end{align}
\end{cor}
\begin{proof}
These follow from the estimates of $b_{i,k}$ and $c_{i,k}$, as $\wpq = \sum_i \wpq[q+1,i] = \sum_{k,i} \nabla \Phi_i^{-1} b_{i,k} \ee^{\ii\lambda_{q+1} k\cdot\Phi_i}$, $\wcq =  \sum_i \wcq[q+1,i]  = \sum_{k,i} c_{i,k}\ee^{\ii\lambda_{q+1} k\cdot\Phi_i}$, and $\nabla\Phi_i\sim I_{3\times 3}$. Specifically, by the disjoint support of $\wpq[q+1,i]$ (see Figure \ref{f:all-eta}),
\[ \|\wpq\|_0 \le \sum_{k\neq0} \sup_{i\ge0} \| \nabla\Phi_i^{-1} b_{i,k}\|_0 \overset{\eqref{e:nabla-phi-i-minus-I3x3},\eqref{e:b-i-k-C0estimate}} \lesssim 2\overline M\delta_{q+1}^{1/2}\sum_{k\neq0} \frac1{|k|^6} \le \frac{M}{32}\delta_{q+1}^{1/2}. \]
For the $C^1$ norm, since $|\nabla \ee^{\ii \lambda_{q+1} k\cdot\Phi_i}|\le 2 \lambda_{q+1} |k|$ by \eqref{e:nabla-phi-i-minus-I3x3},
\begin{align*}
     \|\wpq\|_1 &\le  \frac{M}{32}\delta_{q+1}^{1/2} + \sum_{k\neq0} \sup_{i\ge0} \| \nabla\Phi_i^{-1} b_{i,k} \ee^{\ii \lambda_{q+1}k\cdot \Phi_i} \|_1
     \\
     &\le  \frac{M}{32}\delta_{q+1}^{1/2} + \sum_{k\neq 0}  \ell_q^{-1} 2\sup_{i\ge0}\|b_{i,k}\|_0 + 2\sup_{i\ge0}\| b_{i,k}\|_1 + 4\lambda_{q+1} |k|\sup_{i\ge0}\|b_{i,k} \|_0
     \\
     & \le   \frac{M}{32}\delta_{q+1}^{1/2} +  \sum_{k\neq 0}\frac{(2C_0\ell_q^{-1}+2C_1)\delta_{q+1}^{1/2}}{\ell_q|k|^6} + \frac{4\lambda_{q+1} \overline{M} \delta_{q+1}^{1/2}}{|k|^4}
     \\
     & \le M\lambda_{q+1}  \delta_{q+1}^{1/2}  \left(\frac{1}{32 \lambda_{q+1}} + \frac{(2C_0+2C_1)\sum_{k\neq0}|k|^{-6} }{4\ell_q\lambda_{q+1}} + \frac 1{16} \right),
\end{align*}
where $C_0,C_1$ are the implicit constants coming from the $N=0,1$ cases of \eqref{e:bik-CN}. Note that $\ell_q \lambda_{q+1}\ge \frac1{50}\ell_q \lambda_q^b\ge \frac1{100} \lambda_q^{(1-\beta)(b-1) - \frac{3\alpha}2}$ by the estimates in Section \ref{ss:params}. Hence, \eqref{e:wpq-decay} will follow by choosing $\alpha$ so that
\begin{align}
    \frac{3\alpha}2 < \frac{(1-\beta)(b-1)}2 \label{e:params6-alpha}
\end{align}
and then taking $a\gg1$ to control constants, for instance one could take $a>1+(1600+(C_0+C_1)\sum_{k\neq0}|k|^{-6})^{\frac 2{(1-\beta)(b-1)}}$.

The same choices of constants easily implies the estimate \eqref{e:wcq-decay}, and  the inequality $\lambda_{q+1}\ell_q>4$, from which \eqref{e:wq+1-decay} follows.
\end{proof}

\subsection{Derivation of new stress error term}
By defining $v_{q+1}=\vv_q+w_{q+1}$, we have already implicitly defined the  error term $\div \RR_{q+1}$ through the equation \eqref{e:subsol-euler} with $q$ replaced by $q+1$:
\begin{align*}
    \MoveEqLeft \del_t v_{q+1} + \div(v_{q+1} \otimes v_{q+1})
     \\& = \del_t \vv_q + \div (\vv_q \otimes \vv_q)  + \nabla \pp_q  - \div\RRR_q \\
    & \qquad + \div \RRR_q + \div(w_{q+1}\otimes w_{q+1}) - \nabla \pp_q
    \\
    & \qquad + \del_t w_{q+1} + \vv_q \cdot \nabla w_{q+1} + \div(w_{q+1} \otimes \vv_q).
\end{align*}
The first line of the right-hand side vanishes by \eqref{e:subsol-glued-euler};
for the other terms we group as follows:
\begin{align*}
p_{q+1}(x,t)&\coloneq \pp_q(x,t) - 3\sum_{i=0}^{i_{\textup{max}}} \rho_{q,i}(x,t) +\rho_q(t),
\\
    \Rosc &\coloneq  \mathcal R(  \div (w_{q+1}\otimes w_{q+1})+\div\RRR_q -\nabla p_{q+1}) 
\\
\Rtransport &\coloneq \mathcal R( \del_t w_{q+1} + \vv_q \cdot \nabla w_{q+1}) = \mathcal R(\matd{\vv_q} w_{q+1})
\\
\Rnash &\coloneq \mathcal R \div (w_{q+1}  \otimes \vv_q )= \mathcal R ( w_{q+1}\cdot\nabla \vv_q ) 
\end{align*}
The term $3\sum_i \rho_{q,i}$ comes from the divergence of $\RR_{q,i}$ (see \eqref{e:pressure1}), and $\rho_q$ is subtracted so that $p_{q+1}$ is mean-free; it does not depend on $x$ so disppears on taking the gradient.
Then we set
\[ \RR_{q+1} \coloneq  \Rosc + \Rtransport + \Rnash.\]

\section{Estimates of the stress error terms}
In this section, we show how to achieve the estimate
\[ \|\RR_{q+1}\|_{\alpha} \lesssim \frac{\delta_{q+1}^{1/2}\delta_q^{1/2} \lambda_q}{\lambda_{q+1}^{1-4\alpha}}. \]
The calculations that we omit can be found in \cite[Subsection 6.1]{zbMATH07038033}.
We claim that we can choose parameters such that
\begin{align}
    \frac{\delta_{q+1}^{1/2} \delta_q^{1/2}\lambda_q}{\lambda_{q+1}} \lesssim \frac{\delta_{q+2}}{\lambda_{q+1}^{8\alpha}}. \label{e:estimate-from-parameters0}
\end{align}
which therefore implies $\|\RR_{q+1}\|_{\alpha} \lesssim   \delta_{q+2}\lambda_{q+1}^{-4\alpha}.$ Since we have one extra copy of $\lambda_q^{-\alpha}\ll1$,  $a\gg 1$  gives \eqref{e:RR_q-C0} from \eqref{e:estimate-from-parameters0}.

Since \eqref{e:estimate-from-parameters0} follows if we have $2\beta(b-1)^2-(1-3\beta)(b-1)+8\alpha b < 0$, it suffices to take
\begin{align}
    0<b-1<\frac{1-3\beta}{2\beta}  , \qquad 8\alpha  < -\frac{2\beta(b-1)^2-(1-3\beta)(b-1)}{2b}, \label{e:params5-b-alpha}
\end{align} so that
$
     2\beta(b-1)^2-(1-3\beta)(b-1) +8\alpha b < \frac{2\beta(b-1)^2-(1-3\beta)(b-1)}2 < 0.
$

We will also need the estimate
\begin{align}
    \frac{1}{\lambda_{q+1}^{N-\alpha} \ell_q^{N+\alpha}} &\leq \frac{1}{\lambda_{q+1}^{1-\alpha}} \label{e:estimate-from-parameters}\end{align}
which holds for sufficiently small $\alpha$,  $N=N(b,\beta)$ sufficiently large and large $a$,  independently of $q$ as follows. it is equivalent to the estimate $\lambda_{q+1}^{-N+1} \ell_q^{-N-\alpha} \le 1$. Since $\ell_q^{-\alpha} \le \lambda_q^{2\alpha}$, it suffices to show the stronger estimate $\lambda_{q+1}^{-N+1} \ell_q^{-N}\lambda_q^{2\alpha}\le 1$. Expanding the definition of $\ell_q$ gives
\[ \lambda_{q+1}^{-N(1-\beta)+1}\lambda_q^{N(1-\beta+3\alpha/2) + 2\alpha} \le 1.\]
Since  the exponent of $\lambda_{q+1}$ is negative for $N>2>\frac1{1-\beta}$, and $\lambda_{q+1}\gtrsim \lambda_q^b$, it is enough to enforce 
\[ b - N(1-\beta)(b-1)+1<0 \iff N > \frac{b+2}{(1-\beta)(b-1)}.\] 
This means that in the entire iteration, we only ever use $N$th derivative estimates for a fixed $q$-independent $N$.

\label{s:est-stress}

\begin{prop}[Estimate for $\Rtransport$]
    \[ \|\Rtransport\|_\alpha  \lesssim \frac{\delta_{q+1}^{1/2}\delta_q^{1/2}\lambda_q}{\lambda_{q+1}^{1-4\alpha}}. \]
\end{prop}
\begin{proof}
We separately estimate $\matd{\vv_q}\wpq$ and $\matd{\vv_q}\wcq$. For the first term, we use the fact that $U_{ik} \coloneq (\nabla\Phi_i)^{-1} \rho_{q,i}^{1/2} a_k(\widetilde R_{q,i})\ee^{\ii \lambda_{q+1} k\cdot\Phi_i}$ is Lie-advected up to an error term,
\begin{align*}
    \matd{\vv_q} U_{ik} = (U_{ik}\cdot\nabla) \vv_q + (\nabla\Phi_i)^{-1} (\matd{\vv_q} (\rho_{q,i}^{1/2} a_k(\widetilde R_{q,i})))\ee^{\ii \lambda_{q+1} k\cdot\Phi_i}.
\end{align*}

Importantly, the term $\matd{\vv_q}\ee^{\ii\lambda_{q+1} k\cdot\Phi_i}=0$ does not contribute a bad power of $\lambda_{q+1}$. As $\wpq=\sum_{i,k} U_{ik}$,
\begin{align*}
    \matd{\vv_q}\wpq &=\sum_{i\ge0,k\neq0} ((\nabla\Phi_i)^{-1} \rho_{q,i}^{1/2} a_k(\widetilde R_{q,i})\cdot\nabla)\vv_q \ee^{\ii \lambda_{q+1} k\cdot\Phi_i}
    \\ &\quad + \sum_{i\ge0,k\neq0} (\nabla\Phi_i)^{-1} (\matd{\vv_q} (\rho_{q,i}^{1/2} a_k(\widetilde R_{q,i})))\ee^{\ii \lambda_{q+1} k\cdot\Phi_i}.
\end{align*}
Both these sums are treated similarly. For the first, by Lemma \ref{l:non-stationary-phase}, \eqref{e:nabla-phi-i-CN}, \eqref{e:rho_qi-CN-bd}, \eqref{e:vv_q-bound}, and \eqref{e:mikado-a_k and c_k estimates}:
\begin{align*}
\begin{multlined}
    \|\mathcal R \left(((\nabla\Phi_i)^{-1} \rho_{q,i}^{1/2} a_k(\widetilde R_{q,i})\cdot\nabla)\vv_q \ee^{\ii \lambda_{q+1} k\cdot\Phi_i}
\right) \|_\alpha \\ \lesssim \frac{\lambda_{q} \delta_{q+1}^{1 / 2} \delta_{q}^{1 / 2}}{\lambda_{q+1}^{1-\alpha}|k|^{6}}
    +\frac{\lambda_{q} \delta_{q+1}^{1 / 2} \delta_{q}^{1 / 2}}{\lambda_{q+1}^{N-\alpha} \ell_q^{N+1+3 \alpha}|k|^{6}}
\overset{\eqref{e:estimate-from-parameters}}\lesssim  \frac{\lambda_{q} \delta_{q+1}^{1 / 2} \delta_{q}^{1 / 2}}{\lambda_{q+1}^{1-\alpha}|k|^{6}}.
\end{multlined}
\end{align*}
For the second sum, by Lemma \ref{l:non-stationary-phase}. \eqref{e:nabla-phi-i-CN}, \eqref{e:rho_qi-ddt-CN}, and \eqref{e:mikado-a_k and c_k estimates},
\begin{align*}
\begin{multlined}
     \|\mathcal R \left((\nabla\Phi_i)^{-1} (\matd{\vv_q} (\rho_{q,i}^{1/2} a_k(\widetilde R_{q,i})))\ee^{\ii \lambda_{q+1} k\cdot\Phi_i} \right) \|_{\alpha}
     \\ \lesssim \frac{\delta_{q+1}^{1 / 2}}{\tau_{q} \lambda_{q+1}^{1-\alpha}|k|^{6}}
=\frac{\delta_{q+1}^{1 / 2} \delta_{q}^{1 / 2} \lambda_{q}}{\lambda_{q+1}^{1-\alpha}|k|^{6}} \ell_q^{-2 \alpha} \overset{\eqref{e:estimate-from-parameters}}\lesssim  \frac{\delta_{q+1}^{1 / 2} \delta_{q}^{1 / 2} \lambda_{q}}{\lambda_{q+1}^{1-4 \alpha}|k|^{6}}.
\end{multlined}
\end{align*}
For the last term $\matd{\vv_q} \wcq$ we instead use the estimate \eqref{e:cik-matd-CN} on $\matd{\vv_q}c_{i,k}$ since $\wcq[q+1,i]=\sum_{k\neq0} c_{i,k}\ee^{\ii \lambda_{q+1}k\cdot\Phi_i}$:
\begin{align*}
    &\|\mathcal R (\matd{\vv_q} c_{i,k} )\ee^{\ii \lambda_{q+1}k\cdot\Phi_i} \|_\alpha
    \\
        &\lesssim  \frac{\left\|\matd{\vv_q} c_{i, k}\right\|_{0}}{\lambda_{q+1}^{1-\alpha}}+\frac{\left\|\matd{\vv_q} c_{i, k}\right\|_{N+\alpha}+\left\|\matd{\vv_q} c_{i, k}\right\|_{0}\left\|\Phi_{i}\right\|_{N+\alpha}}{\lambda_{q+1}^{N-\alpha}}
    \\
    &\lesssim \frac{\delta_q^{1/2}|k|^{-6}\ell_q^{-1}\lambda_{q+1}^{-1}\tau_q^{-1}}{\lambda_{q+1}^{1-\alpha}}+\frac{\delta_q^{1/2}|k|^{-6}\ell_q^{-N-\alpha-1}\lambda_{q+1}^{-1}\tau_q^{-1}}{\lambda_{q+1}^{N-\alpha}}
    \\
    &\overset{\eqref{e:estimate-from-parameters}}\lesssim \frac{\delta_{q+1}^{1 / 2}}{\tau_{q} \ell_q \lambda_{q+1}^{2-\alpha}|k|^{6}}
    \lesssim \frac{\delta_{q+1}^{1 / 2}}{\tau_{q} \lambda_{q+1}^{1-\alpha}|k|^{6}}
    \lesssim \frac{\delta_{q+1}^{1 / 2} \delta_{q}^{1 / 2} \lambda_{q}}{\lambda_{q+1}^{1-3 \alpha}|k|^{6}}.
\end{align*}
The required estimate follows by summing in $k$ and the disjointness of the supports of $\wcq[q+1,i]$.
\end{proof}
\begin{prop}[Estimate for $\Rnash$]
    \[ \|\Rnash\|_\alpha  \lesssim \frac{\delta_{q+1}^{1/2}\delta_q^{1/2}\lambda_q}{\lambda_{q+1}^{1-\alpha}}. \]
\end{prop}
\begin{proof}
We write $w_{q+1} =\wpq+\wcq$ and expand $\Rnash$ with \eqref{e:wpqi-bik} and \eqref{e:wcqi-cik}:
\[
\Rnash =\mathcal R \left( \sum_{\substack{k\neq 0\\i\ge0}} (((\nabla \Phi_i)^{-1}b_{i,k})\cdot\nabla \vv_q)  \ee ^{\ii \lambda_{q+1} k\cdot\Phi_i} +  \sum_{\substack{k\neq 0\\i\ge0}} (c_{i,k}\cdot\nabla \vv_q)  \ee ^{\ii \lambda_{q+1} k\cdot\Phi_i} \right).
\]
These sums can be estimated by using Lemma \ref{l:non-stationary-phase}, \eqref{e:vv_q-bound}, \eqref{e:nabla-phi-i-CN}, \eqref{e:bik-CN} and \eqref{e:cik-CN}:%
\begin{align*}
   & \|\mathcal R(((\nabla \Phi_i)^{-1}b_{i,k})\cdot\nabla \vv_q)  \ee ^{\ii \lambda_{q+1} k\cdot\Phi_i} \|_\alpha \\ &  \lesssim \frac{\|\nabla \Phi_{i}^{-1} b_{i, k} \cdot \nabla \bar{v}_{q}\|_{0}}{\lambda_{q+1}^{1-\alpha}}
   +\frac{\|\nabla \Phi_{i}^{-1} b_{i, k} \cdot \nabla \bar{v}_{q}\|_{N+\alpha}\!+\|\nabla \Phi_{i}^{-1} b_{i, k} \cdot \nabla \bar{v}_{q}\|_{0}\|\Phi_{i}\|_{N+\alpha}}{\lambda_{q+1}^{N-\alpha}} \\ & \lesssim  \frac{\lambda_{q} \delta_{q+1}^{1 / 2} \delta_{q}^{1 / 2}}{\lambda_{q+1}^{1-\alpha}|k|^{6}}+\frac{\lambda_{q} \delta_{q+1}^{1 / 2} \delta_{q}^{1 / 2}}{\lambda_{q+1}^{N-\alpha} \ell_q^{N+\alpha}|k|^{6}} \overset{\eqref{e:estimate-from-parameters}}\lesssim \frac{\delta_{q+1}^{1 / 2} \delta_{q}^{1 / 2} \lambda_{q}}{\lambda_{q+1}^{1-\alpha}|k|^{6}}, \text{ and} \\
 &  \| \mathcal{R} ((c_{i, k}\cdot \nabla \bar{v}_{q})\ee^{\ii \lambda_{q+1} k \cdot \Phi_{i}} ) \|_{\alpha} \\ & \lesssim \frac{\left\|c_{i, k} \cdot \nabla \bar{v}_{q}\right\|_{0}}{\lambda_{q+1}^{1-\alpha}}
 +\frac{\left\|c_{i, k} \cdot \nabla \bar{v}_{q}\right\|_{N+\alpha}+\left\|c_{i, k} \cdot \nabla \bar{v}_{q}\right\|_{0}\left\|\Phi_{i}\right\|_{N+\alpha}}{\lambda_{q+1}^{N-\alpha}} \\ & \lesssim \frac{\delta_{q+1}^{1 / 2} \delta_{q}^{1 / 2} \lambda_{q}}{\ell_q \lambda_{q+1}^{2-\alpha}|k|^{6}}+\frac{\delta_{q+1}^{1 / 2} \delta_{q}^{1 / 2} \lambda_{q}}{\ell_q^{N+1-\alpha} \lambda_{q+1}^{N+1-\alpha}|k|^{6}} \overset{\eqref{e:estimate-from-parameters}}\lesssim \frac{\delta_{q+1}^{1 / 2} \delta_{q}^{1 / 2} \lambda_{q}}{\lambda_{q+1}^{1-\alpha}|k|^{6}}.
\end{align*}
Summing in $k$ gives the claimed inequality.
Details can be found in \cite[Subsection 6.1]{zbMATH07038033}.
\end{proof}

\begin{prop}[Estimate for $\Rosc$]
\[ \|\Rosc\|_\alpha  \lesssim \frac{\delta_{q+1}^{1/2}\delta_q^{1/2}\lambda_q}{\lambda_{q+1}^{1-3\alpha}}. \]
\end{prop}
\begin{proof}
Note that
$w_{q+1}\otimes w_{q+1}=\wpq\otimes\wpq + (\wpq \otimes \wcq+\wcq \otimes \wpq+\wcq \otimes \wcq)$. We will momentarily show
\begin{align}
    \left\|\mathcal R\div\left( \wpq\otimes\wpq - \RRR_q -3\sum_{i=0}^{i_{\textup{max}}} \nabla \rho_{q,i}  \right)\right\|_{\alpha} \lesssim \frac{\delta_{q+1}^{1/2}\delta_q^{1/2}\lambda_q}{\lambda_{q+1}^{1-3\alpha}} \label{e:Rosc-main}.
\end{align}
The other terms are sufficiently small because $\wcq$ is small: interpolating the estimates \eqref{e:wpq-decay} and \eqref{e:wcq-decay} gives
$ \|\wpq\|_\alpha \lesssim \delta_{q+1}^{1/2}\lambda_{q+1}^{\alpha}$, $\|\wcq\|_\alpha \lesssim \delta_{q+1}^{1/2}\ell_q^{-1}\lambda_{q+1}^{-1+\alpha}.$
Therefore,
\begin{align*}
    \MoveEqLeft\| w_{q+1}\otimes w_{q+1}-\wpq\otimes\wpq\|_\alpha
    \\
    &\lesssim \|\wpq\|_0\|\wcq\|_\alpha + \|\wpq\|_\alpha\|\wcq\|_0 + \|\wcq\|_\alpha\|\wcq\|_0
    \\
    & \le \delta_{q+1}\ell_q^{-1}\lambda_{q+1}^{-1+\alpha} + \delta_{q+1}^{1/2}\ell_q^{-2}\lambda_{q+1}^{-2+\alpha}
    \\
    &\lesssim \frac{\delta_{q+1}^{1/2} \delta_q^{1/2}\lambda_q^{1+3\alpha/2} }{\lambda_{q+1}^{1-\alpha}}
    \le \frac{\delta_{q+1}^{1/2}\delta_q^{1/2}\lambda_q}{\lambda_{q+1}^{1-3\alpha}}.
\end{align*}
To show \eqref{e:Rosc-main}, we expand out $\wpq\otimes\wpq=\sum_{ij} \wpq[q+1,i]\otimes \wpq[q+1,j]$. Since $\eta_i$ have disjoint supports, it follows that $\wpq[q+1,i]\otimes \wpq[q+1,j]=0$ if $i\neq j$. Hence
\begin{align}
\wpq[q+1,i]\otimes \wpq[q+1,i] &= \wpq[q+1,i](\wpq[q+1,i])^\TT   \notag
\\&=  \rho_{q,i} \nabla \Phi_i^{-1} (W\otimes W)(\widetilde R_{q,i},\lambda_{q+1}\Phi_i ) \nabla \Phi_i^{-\TT} \notag
    \\& \overset{\mathclap{\eqref{e:mikado-C_k}}}= \rho_{q,i}\nabla \Phi_i^{-1} \widetilde R_{q,i} \nabla \Phi_i^{-\TT}  + \sum_{k\neq0} \rho_{q,i} \nabla \Phi_i^{-1}C_k(\widetilde R_{q,i})\nabla \Phi_i^{-\TT}
        \ee^{\ii \lambda_{q+1}k\cdot \Phi_i} \notag
    \\& \overset{\mathclap{\eqref{e:defn-Rqi}}}= 
    \rho_{q,i}I_{3\times 3} -  \RRR_q  + \sum_{k\neq0}  \rho_{q,i}  \nabla \Phi_i^{-1}C_k(\widetilde R_{q,i})\nabla \Phi_i^{-\TT}
        \ee^{\ii \lambda_{q+1}k\cdot \Phi_i} .\label{e:wpq-otimes-wpq-expanded}
\end{align}
We stress that the above step works despite our modified definition of $\rho_{q,i}$. We need to compute $\div(\wpq[q+1,i]\otimes \wpq[q+1,i])$. Observe that $\div (Mf)=(\div M)f + M\nabla f$, valid for matrix valued $M$ and scalar $f$, and the chain rule $\nabla (g(k\cdot\Phi_i))=\nabla \Phi^\TT k g'$, when $g:\mathbb R\to\mathbb R$. Putting $f=\ee^{\ii\lambda_{q+1}k\cdot \Phi_i}$, $g=\exp,$ we see that when the derivative falls on $\ee^{\ii \lambda_{q+1}k\cdot \Phi_i}$, this left factor of $\nabla\Phi^\TT$ cancels with the right factor of $\nabla\Phi_i^{-\TT}$, allowing the use of $C_k\cdot k=0$ in \eqref{e:cancellation}:
\begin{align}
 \MoveEqLeft \div(\wpq[q+1,i]\otimes \wpq[q+1,i]+\div\RRR_q)\notag  \\ &= 3\nabla \rho_{q,i} + \sum_{i\ge 0, k\neq0} \div\Big(\rho_{q,i}\nabla \Phi_i^{-1}C_k(\widetilde R_{q,i})\nabla \Phi_i^{-\TT}
        \Big)\ee^{\ii \lambda_{q+1}k\cdot \Phi_i} \label{e:pressure1}
\\& +\lambda_{q+1} \sum_{i\ge 0, k\neq0} \rho_{q,i} \nabla \Phi_i^{-1}\underbrace{C_k(\widetilde R_{q,i})\nabla \Phi_i^{-\TT} \nabla \Phi_i^{\TT} k}_{=C_k(\widetilde R_{q,i})k\equiv 0}
        \ee^{\ii \lambda_{q+1}k\cdot \Phi_i}
 \label{e:cancellation}	
\end{align}
The first term of \eqref{e:pressure1} is specifically subtracted in \eqref{e:Rosc-main}. To finish the calculation, we set $a\coloneq \div( \rho_{q,i}\nabla \Phi_i^{-1}C_k(\widetilde R_{q,i})\nabla \Phi_i^{-\TT})$ and apply Lemma \ref{l:non-stationary-phase} to find \eqref{e:Rosc-main}, as needed:
\begin{align*}
        &\left\|\mathcal R\div\left( \wpq\otimes\wpq - \RRR_q -3\sum_{i=0}^{i_{\textup{max}}} \nabla \rho_{q,i}  \right)\right\|_{\alpha}
        \\
        &=
        \sum_{i\ge0,k\neq 0} \frac{\|a\|_0}{\lambda_{q+1}^{1-\alpha}} + \frac{\|a\|_{N+\alpha}}{\lambda_{q+1}^{N-\alpha}} + \frac{\|a\|_0\|\Phi_i\|_{N+\alpha}}{\lambda_{q+1}^{N-\alpha}}
        \\
        &\lesssim \sum_{i\ge0,k\neq0} \frac{\delta_{q+1}}{\ell_q\lambda_{q+1}^{1-\alpha}|k|^6}
        \lesssim \frac{\delta_{q+1}^{1/2}\delta_q^{1/2}\lambda_q}{\lambda_{q+1}^{1-3\alpha}} .\qedhere \end{align*}
 \end{proof}
\section{Energy iteration}
\label{s:energy-iteration}
\begin{prop}[Energy estimate for $v_{q+1}$]
\label{p:energy}For $t\in[1-\tau_q,T]$,
    \begin{align}
    \bigg|e(t)-\int_{\mathbb T^3} |v_{q+1}|^2\dd x -\frac{\delta_{q+2}}2\bigg| \lesssim \frac{\delta_q^{1/2}\delta_{q+1}^{1/2} \lambda_q^{1+2\alpha}}{\lambda_{q+1}}\le\delta_{q+2} . \label{e:energy-estimate}
\end{align}

\end{prop}
\begin{proof}
Notice that we only want to control the energy for times $t\ge1-\tau_q$. For such times, the $\zeta$ term (defined in Section \ref{ss:energy-gap}) in the energy gap $\rho_{q,i}$ is identically zero, i.e. we can write \eqref{e:energy-q-i-gap} as
\[
 \rho_{q,i} (x,t) = \frac{\eta^2_i(x,t)}{ \sum_{i= 0}^{i_{\textup{max}}} \int_{\mathbb T^3}\eta_i (\tilde x,t)\dd \tilde x }\rho_q(t).
\]
Hence, for $t\in[1-\tau_q,T]$, our $\rho_{q,i}$ matches the formula in \cite[Section 5.2]{zbMATH07038033}, and therefore we can simply use their proof, which we paraphrase here for the reader's convenience.

We decompose the total energy into three parts,
\[\int_{\mathbb T^3} |v_{q+1}|^2\dd x = \int_{\mathbb T^3} |\vv_q|^2\dd x  + 2 \int_{\mathbb T^3} w_{q+1}\cdot \vv_q \dd x +\int_{\mathbb T^3} |w_{q+1}|^2\dd x \eqcolon I_1+I_2+I_3.\]
Note that $I_2$ is sufficently small by integration by parts, since $w_{q+1}$ is a curl, and we have the estimates \eqref{e:nabla-phi-i-CN}, \eqref{e:bik-CN}, and \eqref{e:vv_q-bound}:
\[ \left|\int_{\mathbb T^3} \wpq\cdot \vv_q \dd x \right| \lesssim \frac1{\lambda_{q+1}} \sum_{i\ge0}\sum_{k\neq0} \left\|\nabla \Phi_i^\TT \frac{\ii k \times b_k}{|k|^2} \right\|_0 \| \vv_q \|_1 \lesssim \frac{\delta_q^{1/2} \delta_{q+1}^{1/2} \lambda_q}{\lambda_{q+1}}.  \] Now we will show that  $I_3\approx\int_{\mathbb T^3}|\wpq|^2\dd x \approx 3\rho_q$, which was designed in \eqref{e:energy-gap} specifically to cancel with $I_1=\int_{\mathbb T^3} |\vv_q|^2\dd x $, leaving behind an error term of size $\approx\delta_{q+2}$.

    To estimate $I_3$, we rewrite
    \[ I_3 = \int_{\mathbb T^3}|\wpq|^2 \dd x + \int_{\mathbb T^3} \wcq\cdot(\wpq+\wcq)\dd x\eqcolon I_{31}+I_{32}.\] $I_{32}$ is an error term, from the estimates \eqref{e:wcq-decay} and \eqref{e:wq+1-decay} and the definition \eqref{e:ell} of $\ell_q$:
    \[ I_{32} \lesssim \frac{\delta_{q+1}}{\ell_q \lambda_{q+1}} = \frac{\delta_q^{1/2}\delta_{q+1}^{1/2} \lambda_q^{1+2\alpha}}{\lambda_{q+1}}. \] For $I_{31}$, we take the trace of \eqref{e:wpq-otimes-wpq-expanded} and use $\tr\RRR_q=0$ to find
\begin{align*}
    &\int_{\mathbb T^3} |\wpq|^2 \dd x = \int_{\mathbb T^3} \tr (\wpq\otimes \wpq) \dd x
    \\
    &= 3\rho_{q} + \tr \Big( \int_{\mathbb T^3} \sum_{i\ge0,k\neq0} \rho_{q,i} \nabla \Phi_i^{-1} C_k(\widetilde R_{q,i}) \nabla \Phi_i^{-\TT } \ee ^{\ii \lambda_{q+1} k\cdot\Phi_i} \dd x\Big),
\end{align*}
since $\rho_q = \int_{\mathbb T^3} \sum_{i\ge0} \rho_{q,i}\dd x$.  From \eqref{e:rho_qi-C0-bd}, \eqref{e:tilde-Rqi-N}, \eqref{e:cik-CN}, and \eqref{e:nabla-phi-i-CN}, we see that
\[ \|\nabla \Phi_i^{-1} C_k(\widetilde R_{q,i}) \nabla \Phi_i^{-\TT }\|_N\lesssim \delta_{q+1}\ell_q^{-N}. \]
Since $\ell_q\ll \lambda_{q+1}$, repeated integration by parts (formalised by Lemma \ref{l:non-stationary-phase}) shows that this term is arbitrarily small:
\[\left|\int_{\mathbb T^3} |\wpq|^2 \dd x - 3\rho_q \right| \lesssim \sum_{k\neq0} \frac{\delta_{q+1}\ell_q^{-N}}{\lambda_{q+1}^N|k|^N} \]
which is finite for $N\ge 4$.  Taking $N=N(b,\beta)\gg 1$ as in Section \ref{s:est-stress} so that 
\[ \frac{\delta_{q+1}\ell_q^{-N}}{\lambda_{q+1}^N} \le \frac{\delta_{q+1}\delta_q^{1/2}\lambda_q}{\lambda_{q+1}}\]
gives the required bound \eqref{e:energy-estimate}, since $3\rho_q = e(t)-\frac{\delta_{q+2}}2-I_1$ by \eqref{e:energy-gap}.
\end{proof}
\appendix
\section{\Holder spaces and estimates} \label{s:holder}
For $N\in\mathbb Z_{\ge0}$, $C^{N}(X)$ denotes the space of $N$ times differentiable functions with the norm
$ \|f\|_{C^{N}(X)} \coloneq \|f\|_{L^\infty(X)} + \sum_{|\sigma|\le N} \| D^\sigma f\|_{L^\infty(X)}.$
For $N\in\mathbb Z_{\ge0}$ and $\alpha\in(0,1)$, $C^{N+\alpha}(X)$ denotes the subspace of $C^N(X)$ whose $N$th derivatives are $\alpha$-\Holder continuous, with the norm
\[ \|f\|_{C^{N+\alpha(X) }} \coloneq \|f\|_{C^N} + \sum_{|\sigma|=N} [D^\sigma f]_{C^\alpha(X)},\]
 where $[f]_{C^\alpha(X)} \coloneq \sup_{x,y\in X : x\neq y} \frac{|f(x)-f(y)|}{|x-y|^\alpha}$ is the \Holder seminorm.

\begin{prop}
Let $f$ be a smooth function $f:\mathbb T^3\to\mathbb R$, and let the mollifier $\psi_\epsilon$ be as defined in \eqref{e:defn-mollifier-x}. For $\epsilon\ll 1$, $\alpha\in(0,1)$, $\beta\in(0,1)$, $N\ge 0$, we have the estimates
\begin{align}
\|f*\psi_\epsilon \|_{C^{N+\alpha}}&\lesssim \|f\|_{C^\alpha}\epsilon^{-N} \label{e:convolution-holder-estimate}    ,
\\
\|f-f*\psi_\epsilon\|_{C^\alpha} &\lesssim [f]_{C^\beta} \epsilon^{\beta-\alpha}     \label{e:convolution-diff-holder-estimate} \quad \text{(if $\alpha\le \beta$)},
\\
\|f-f*\psi_\epsilon\|_{C^\alpha} &\lesssim [\nabla f]_{C^\beta} \epsilon^{1+\beta-\alpha}    \label{e:convolution-diff-holder-estimate2}.
\end{align}
\end{prop}

\begin{proof}
These estimates are standard; for completeness, we give their short proofs. For \eqref{e:convolution-holder-estimate}, we have for $\epsilon\ll 1$,
\[ \|f*\psi_\epsilon\|_{C^{N+\alpha}} \le \|f\|_{C^{0}} + \|f*\nabla^N\psi_\epsilon\|_{C^{\alpha}} \lesssim \epsilon^{-N}\|f\|_{C^{\alpha}}.\]
For \eqref{e:convolution-diff-holder-estimate}, we first have the $L^\infty$ estimate
\[ |f(x) - f*\psi_\epsilon (x) | \leq \int_{\mathbb B_\epsilon(0)} |f(x) - f(x-y)| \psi_\epsilon (y) \dd y \lesssim [f]_{C^{\beta}} \epsilon^\beta. \]
We also have
$ |f(x) - f*\psi_\epsilon (x) - (f(z) - f*\psi_\epsilon (z)) | \lesssim [f]_{C^\beta} |x-z|^\beta,$
so that $[f-f*\psi_\epsilon]_{C^\beta} \lesssim [f]_{C^{\beta}}$ and hence that $ \|f-f*\psi_\epsilon\|_{C^\beta}\lesssim [f]_{C^\beta}$. So \eqref{e:convolution-diff-holder-estimate} follows by  interpolating
$\|f-f*\psi_\epsilon\|_{C^\alpha} \le \|f-f*\psi_\epsilon\|_{C^0}^{(\beta-\alpha)/\beta} \|f-f*\psi_\epsilon\|_{C^\beta}^{\alpha/\beta}.$

The estimate \eqref{e:convolution-diff-holder-estimate2} similarly follows from the interpolation inequality
\[\|f-f*\psi_\epsilon\|_{C^\alpha} \le \|f-f*\psi_\epsilon\|_{C^0}^{(1+\beta-\alpha)/(1+\beta)} \|f-f*\psi_\epsilon\|_{C^{1+\beta}}^{\alpha/(1+\beta)}\]
  once we prove the corresponding $L^\infty$ bound. For this, we use the fact that symmetry of the mollifier implies that
$  \nabla f(x)\cdot \int_{B_\epsilon(0)} y \psi_\epsilon(y)\dd y=0.$
Hence,
\[ f(x) - f*\psi_\epsilon(x) = \int_{\mathbb B_\epsilon(0)} \Big(f(x) - \nabla f(x) \cdot y - f(x-y) \Big) \psi_\epsilon(y)\dd y. \]
Since the fundamental theorem of calculus implies
$  |f(x+y)-f(x)-\nabla f(x)\cdot y| \leq \int_0^1 |\nabla f(x+ty)-\nabla f(x)| |y| \dd t \leq |y|^{1+\beta}[\nabla f]_{C^\beta} \frac1{1+\beta},$
we obtain
\[ \|f-f*\psi_\epsilon\|_{C^0} \lesssim [\nabla f]_{C^\beta}\epsilon^{1+\beta}.\]
 Since by \eqref{e:convolution-diff-holder-estimate} we also have
$
 \|\nabla f - \nabla f*\psi_\epsilon\|_{C^\beta} \lesssim [\nabla f]_{C^\beta},
$ 
the inequality \eqref{e:convolution-diff-holder-estimate2} follows.
\end{proof}

\section{Inverse Divergence Operator}
\label{s:inverse-div}
We recall  \cite[Def. 1.4]{zbMATH06456007} the following  operator of order $-1$,
\begin{align*}
    \mathcal R u = -(-\Delta)^{-1} (\nabla u + \nabla u^\TT ) - \frac12(-\Delta)^{-2} \nabla^2 \nabla\cdot u  + \frac12 (-\Delta)^{-1} (\nabla\cdot u )I_{3\times 3},
\end{align*}
where $U=(-\Delta)^{-1}v$ is the mean-free solution to $-\Delta U = u-\dashint_{\mathbb T^3} u.$
 It is a matrix-valued right inverse of the divergence operator for mean-free vector fields, in the sense that
\[ \div\mathcal R u= u - \dashint_{\mathbb T^3} u. \]
In addition, $\mathcal Ru$ is traceless and symmetric, satisfies the Schauder estimate
\begin{align}
      \| \mathcal R u\|_{C^{1+\alpha}} \lesssim \|u\|_{C^\alpha},  \label{e:schauder}
\end{align}
and the following non-stationary phase type estimates:
\begin{lem}\label{l:non-stationary-phase} If $a\in C^\infty(\mathbb T^3)$ and $\Phi\in  C^\infty (\mathbb T^3 ; \mathbb R^3)$  such that
\[ 1 \lesssim |\nabla \Phi|  \lesssim 1,\]
then
\[ \left| \int_{\mathbb T^3} a(x) \ee ^{\ii k \cdot \Phi} \dd x \right| \lesssim  \frac{\|a\|_{C^N}+\|a\|_{C^0}\|\Phi\|_{C^N}}{|k|^N},\]
and
\[
    \|\mathcal R( a \ee^{\ii k\cdot \Phi}) \|_{C^\alpha}
    \lesssim  \frac{\|a\|_{C^0}}{|k|^{1-\alpha}} + \frac{\|a\|_{C^{N+\alpha}} + \|a\|_{C^0} \|\Phi\|_{C^{N+\alpha}}} {|k|^{N-\alpha}}.
    \]
\end{lem}

\section*{Acknowledgements}
This work was supported by the National Key Research and Development Program of China (No. 2020YFA0712900) and NSFC Grant 11831004.

\printbibliography

\end{document}